\documentclass[11pt]{article}
\usepackage{pxfonts}
\usepackage{yfonts}
\usepackage{dsfont}
\usepackage{graphicx}
\usepackage{relsize}

\def\bel{\begin{equation}\label}
\def\eeq{\end{equation}}
\def\ds{\displaystyle}
\def\endproof{\hphantom{MM}
\hfill\llap{$\square$}\goodbreak}

\def\mt{\longrightarrow}
\def\v{\vskip 1em}

\def\ve{\varepsilon}
\def\R{\mathbb R}
\def\Z{\mathbb Z}
\def\C{\mathfrak{B}}

\def\N{{\bf N}}

\def\S{{\bf S}}

\def\F{\mathfrak{F}}

\def\O{{\bf O}}
\def\Q{{\bf Q}}
\def\D{{\bf D}}

\def\B{{\bf B}}
\def\H{{\bf H}}
\def\L{{\bf L}}
\def\T{{\bf T}}
\def\K{{\bf K}}
\def\U{\mathcal{U}}
\def\V{{\bf V}}

\def\p{{\partial}}

\def\b{{\bf b}}
\def\i{{\bf i}}
\def\Tilde{\widetilde}
\def\Hat{\widehat}

\def\bar{\overline}
\def\supp{{\bf supp}}

\def\dist{{\bf dist}}

\def\Cup{{\bigcup}}

\def\alpha{\alphaup}
\def\beta{\betaup}
\def\gamma{\gammaup}
\def\delta{\deltaup}
\def\xi{{\xiup}}
\def\eta{{\etaup}}
\def\tau{{\tauup}}
\def\rho{{\rhoup}}
\def\phi{{\phiup}}
\def\psi{{\psiup}}
\def\lambda{{\lambdaup}}
\def\omega{\omegaup}
\def\varphi{{\varphiup}}
\def\gamma{{\gammaup}}
\def\c{{\bf c}}

\def\0{{\bf 0}}

\newtheorem{lemma}{Lemma}[section]

\newtheorem{remark}{Remark}[section]

\begin{document}
 \[\begin{array}{cc}\hbox{\LARGE{\bf Regularity of Fourier integrals on product spaces}}
  \end{array}\]

 \[\hbox{Chaoqiang Tan~~~~~and ~~~~~Zipeng Wang}\]

 \begin{abstract}
We study a family of Fourier integral operators by allowing their symbols to satisfy  a multi-parameter  differential inequality in $\R^\N$. We show that these operators of order $-{\N-1\over 2}$ are bounded from the classical, atom decomposable $\H^1$-Hardy space to $\L^1(\R^\N)$. Consequently, 
we obtain a sharp $\L^p$-regularity result of Seeger, Sogge and Stein.

\end{abstract}
\section{Introduction}
\setcounter{equation}{0}
The theory of Fourier integrals begins from earlier independent researches of Fourier and Cauchy on the propagation of waves. Over the several past decades, it becomes one of the core subjects in harmonic analysis.

A Fourier integral operator $\F$ is defined by
\bel{Ff}
\F f(x)~=~\int_{\R^\N}e^{2\pi\i \Phi(x,\xi)} \sigma(x,\xi) \Hat{f}(\xi)d\xi
\eeq
where
$\sigma(x,\xi)\in\mathcal{C}^\infty(\R^\N\times\R^\N)$ has a compact support in $x$. 

The phase $\Phi$ is  real-valued, homogeneous of degree $1$ in $\xi$ and smooth in $(x,\xi)$ for $\xi\neq0$. Moreover, it satisfies   
 the non-degeneracy condition
\bel{nondegeneracy}
 \det\left[{\p^2\Phi\over \p x\p \xi}\right]\left(x,\xi\right)~\neq~0,\qquad \xi\neq0
 \eeq
on the support of $\sigma(x,\xi)$.

$\diamond$  {\small Throughout,  $\C$ is a generic constant depending on its subindices.}

We say $\sigma\in S^{-\rho}$ if 
\bel{class}
\left|\p_\xi^\alpha\p_x^\beta \sigma(x,\xi)\right|~\leq~\C_{\alpha~\beta}  \left({1\over 1+|\xi|}\right)^{|\alpha|+\rho}
\eeq
for every multi-indices $\alpha,\beta$.

For $\sigma\in S^0$, 
$\F$ defined in (\ref{Ff})-(\ref{nondegeneracy}) is  bounded on $\L^2(\R^\N)$ as shown by  Eskin \cite{Eskin} and H\"{o}rmander \cite{Hormander}. On the other hand,  it is well known that $\F$ of order zero is not bounded on $\L^p(\R^\N)$, $p\neq2$. A natural question asks for what is the optimal $\L^p$-regularity of $\F$ when $\sigma$ having a negative order.
The regarding $\L^p$-estimate was first investigated by   Duistermaat and H\"{o}rmander  \cite{Duistermaat-Hormander}  and then by  Colin de Verdi\'{e}re and Frisch \cite{Colin-Frisch},  Brenner \cite{Brenner}, Peral \cite{Peral}, Miyachi \cite{Miyachi}, Beals \cite{Beals} and eventually proved by Seeger, Sogge and Stein \cite{S.S.S}.

{\bf Theorem One:~~Seeger, Sogge and Stein (1991)}\\ 
 {\it Let $\F$ defined in (\ref{Ff})-(\ref{nondegeneracy}). Suppose $\sigma\in S^{-\rho}$ for $0\leq\rho\leq(\N-1)/2$. We have
  \bel{Result ONE H}
 \left\| \F f\right\|_{\L^1(\R^\N)} ~\leq~\C~\left\| f\right\|_{\H^1(\R^\N)},\qquad 
  \sigma\in S^{-{\N-1\over 2}}
 \eeq 
 and
  \bel{Result ONE}
  \begin{array}{cc}\ds
  \left\| \F f\right\|_{\L^p(\R^\N)}~\leq~\C_p~\left\| f\right\|_{\L^p(\R^\N)},\qquad 

\left|  {1\over p}- {1\over 2}\right|~\leq~{\rho\over \N-1}.
  \end{array}
 \eeq}

$\H^1(\R^\N)$ in (\ref{Result ONE H}) is the $\H^1$-Hardy space  investigated by Fefferman and Stein \cite{FC.S}. 
\begin{remark}
The $\L^p$-norm inequality in (\ref{Result ONE}) is sharp. Let $a\in\mathcal{C}^\infty_o(\R^\N)$ and $a(x)\neq0$ for $|x|=1$. Consider
\bel{example}
\sigma(x,\xi)~=~a(x)\left(1+|\xi|\right)^{-\rho},\qquad \Phi(x,\xi)~=~x\cdot\xi+|\xi|.
\eeq
Then, $\F$ defined in (\ref{Ff})-(\ref{nondegeneracy}) is not bounded on $\L^p(\R^\N)$ for  $\left|{1\over p}-{1\over 2}\right|>{\rho\over \N-1}$.
See {\bf 6.13}, chapter IX in the book of Stein \cite{Stein}.
\end{remark}
 
In this paper, we consider the Fourier integral operator $\F$ defined in (\ref{Ff})-(\ref{nondegeneracy}) but  allow $\sigma(x,\xi)$ to satisfy a  multi-parameter differential inequality.

Let $\N=n+m$ and $\xi=(\tau,\lambda)\in\R^{n}\times\R^{m}$. 
We say $\sigma\in\S^{-\rho}$ if
\bel{Class}
\left|\p_\tau^\alpha\p_\lambda^\beta\p_x^\gamma \sigma(x,\xi)\right|~\leq~\C_{\alpha~\beta~\gamma}~ \left({1\over 1+|\tau|}\right)^{|\alpha|+\rho_1} \left({1\over 1+|\lambda|}\right)^{|\beta|+\rho_2}
\eeq
for every multi-indices $\alpha, \beta$ and $\gamma$ where
\bel{rho_1, rho_2}
\begin{array}{cc}\ds
\rho_1~=~\rho_2~=~0\qquad \hbox{if}\qquad \rho~=~0;
\\\\ \ds
\rho_1~>~\left[{n-1\over \N-1}\right]\rho,\qquad \rho_2~>~\left[{m-1\over \N-1}\right]\rho,\qquad \rho_1+\rho_2~=~\rho\qquad\hbox{if}\qquad \rho~>~0.
\end{array}
\eeq
The study of certain operators  that  commute with a  family of multi-parameter dilations,  dates back to the time of  Jessen, Marcinkiewicz and Zygmund.  
The $\L^p$-theory regarding to the strong maximal function and singular integrals defined on product spaces have been established, for example  by Cordoba and Fefferman \cite{Cordoba-Fefferman}, Fefferman and Stein \cite{R-F.S}, Journ\'{e} \cite{Journe'} and 
M\"{u}ller, Ricci and Stein \cite{M.R.S}. 

In order to explore the possible extension  on $0<p\leq1$, a product version of Hardy spaces, denoted by $\H^p\times \H^p(\R^n\times\R^m)$ was introduced by  M.~Malliavin and P.~Malliavin \cite{Malliavin} and Gundy and Stein \cite{Gundy-Stein}. Such a Hardy space defined by the corresponding multi-parameter Littlewood-Paley inequality  cannot be characterized in terms of  "rectangle atoms". See the counter-example of Carleson \cite{Carleson}. Surprisingly, by testing on any single rectangle atom and using a geometric covering lemma due to Journ\'{e} \cite{Journe'}, Robert Fefferman \cite{R.Fefferman} is able to conclude that a 2-parameter Calder\'{o}n-Zygmund operator defined on $\R^n\times\R^m$ is bounded from $\H^p\times \H^p(\R^n\times\R^m)$ to $\L^p(\R^n\times\R^m)$. The relevant $\H^p\times \H^p(\R^n\times\R^m)\mt\H^p\times \H^p(\R^n\times\R^m)$ theorem is proved more recently by Han, Lee, Lin and Ying \cite{Han}. This area remains largely open to Fourier integral operators.

{\bf Theorem Two}~~
 {\it Let $\F$ defined in (\ref{Ff})-(\ref{nondegeneracy}). Suppose $\sigma\in \S^{-\rho}$ for $0\leq\rho\leq(\N-1)/2$. We have
  \bel{Result A H}
 \left\| \F f\right\|_{\L^1(\R^\N)} ~\leq~\C_{\rho_1~\rho_2}~\left\| f\right\|_{\H^1(\R^\N)},\qquad \sigma\in \S^{-{\N-1\over 2}}
 \eeq 
 and
  \bel{Result A}
  \begin{array}{cc}\ds
  \left\| \F f\right\|_{\L^p(\R^\N)}~\leq~\C_{\rho_1~\rho_2~p}~\left\| f\right\|_{\L^p(\R^\N)},\qquad 
 \left|{1\over p}- {1\over 2}\right|~\leq~{\rho\over \N-1}.
  \end{array}
 \eeq
}
 
\begin{remark}  $\H^1(\R^\N)$ in (\ref{Result A H}) is the same $\H^1$-Hardy space  introduced  by Fefferman and Stein \cite{FC.S}.
Furthermore, it has a characterization of atom decomposition established by Coifman \cite{Coifman}.
\end{remark}
In contrast to the $\H^p\times \H^p(\R^n\times\R^m)\mt\L^p(\R^n\times\R^m)$ regularity estimate  discussed earlier, 
the Fourier integral operator $\F$ with $\sigma\in\S^{-{\N-1\over 2}}$ is an example of multi-parameter operators bounded from the classical, atom decomposable $\H^1(\R^\N)$ to $\L^1(\R^\N)$.

Our proof of {\bf Theorem Two} is developed on a mixture of two different types of dyadic decomposition. In addition to the ingenious construction due to Seeger, Sogge and Stein \cite{S.S.S}, we  further decompose the frequency space  into an infinitely many dyadic cones.  

In the next section, we introduce this new framework consisting of partial operators  defined on every dyadic cone. Moreover, these operators satisfy some $\L^p\mt\L^2$-estimates, stated as {\bf Proposition One}.  
 We prove {\bf Theorem Two} in section 3. Our key result is about the kernels of   regarding partial operators    enjoy certain majorization properties, accumulated into {\bf Proposition Two}.

Section 4 and 5 are devoted to the proof of {\bf Proposition Two}. 
In section 6, we conclude $\F^*\colon\H^1(\R^\N)\mt\L^1(\R^\N)$ for $\sigma\in\S^{-{\N-1\over 2}}$. In section 7, we prove {\bf Proposition One}.

\section{Cone decomposition on $\R^n\times\R^m$}
\setcounter{equation}{0}
Let $\varphi\in\mathcal{C}^\infty_o(\R)$ such that $\varphi(t)=1$ if $|t|\leq1$ and $\varphi(t)=0$ for $|t|>2$.
Define
\bel{delta_t}
\begin{array}{cc}\ds
 \delta_\ell(\xi)~=~\varphi\left[2^{-\ell}{|\tau|\over|\lambda|}\right]-\varphi\left[2^{-\ell+1}{|\tau|\over|\lambda|}\right],\qquad\ell\in\Z
 \end{array}
\eeq
whose support is contained in the dyadic cone
\bel{Cone}
\Lambda_\ell~=~ \Bigg\{  (\tau,\lambda)\in\R^n\times\R^m~\colon~      2^{\ell-1}~<~ {|\tau|\over|\lambda|}~<~2^{\ell+1} \Bigg\}.
\eeq
From  direct computation, we find
\bel{delta_l Diff Ineq}
\left| \p_\tau^\alpha\p_\lambda^\beta   \delta_\ell(\xi)\right|~\leq~\C_{\alpha~\beta}~  |\tau|^{-|\alpha|}|\lambda|^{-|\beta|},\qquad \ell\in\Z
\eeq
for every multi-indices $\alpha,\beta$.
\begin{remark} Given $\ell\in\Z$, $\delta_\ell$ is a Marcinkiewicz multiplier. The regarding convolution operator is bounded on $\L^p(\R^\N)$ for $1<p<\infty$.
\end{remark}

On the other hand, define
\bel{phi_j}
\phi(\xi)~=~\varphi(|\xi|),\qquad \phi_j(\xi)~=~\varphi\left[ 2^{-j} |\xi|\right]-\varphi\left[2^{-j+1}|\xi|\right],\qquad j\in\Z.
\eeq
Let
\bel{F_o}
\begin{array}{cc}\ds
 \F_o f(x)~=~\int_{\R^\N} f(y)\Omega_o (x,y) dy,
\\\\ \ds
\Omega_o(x,y)~=~
\int_{\R^\N} e^{2\pi\i \big[\Phi(x,\xi)-y\cdot\xi \big]}\delta_\ell(\xi)\phi(\xi)\sigma(x,\xi)d\xi.
\end{array}
\eeq
Observe that $\Omega_o(x,y)$ has a $x$-compact support and $\|\Omega_o(x,\cdot)\|_{\L^\infty(\R^\N)}\leq\C$. We clearly have
\bel{F_o L^1}
\left\| \F_o f\right\|_{\L^1(\R^\N)}~\leq~\C\left\| f\right\|_{\L^1(\R^\N)},\qquad \sigma\in\S^0.
\eeq
Given $j>0$ and $\ell\in\Z$, we define
\bel{Delta_lj F}
\begin{array}{cc}\ds
\F_{\ell j} f(x)~=~\int_{\R^\N} f(y)\Omega_{\ell j} (x,y) dy,
\\\\ \ds
\Omega_{\ell j}(x,y)~=~
\int_{\R^\N} e^{2\pi\i \big[\Phi(x,\xi)-y\cdot\xi \big]}\delta_\ell(\xi)\phi_j(\xi)\sigma(x,\xi)d\xi.
\end{array}
\eeq
Furthermore, consider
\bel{delta_j+-}
\begin{array}{lr}\ds
\mathfrak{S}_j^+(\xi)~=~\sum_{\ell\ge j} \delta_\ell(\xi)
~=~ 1-\varphi\left[2^{-j+1}{|\tau|\over|\lambda|}\right],
\\\\ \ds
\mathfrak{S}_j^-(\xi)~=~\sum_{\ell\leq- j} \delta_\ell(\xi)
~=~\varphi\left[2^j{|\tau|\over|\lambda|}\right].
\end{array}
\eeq
Suppose that at least one of $\alpha, \beta$ is non-zero. We find 
$\p_\tau^\alpha\p_\lambda^\beta \left[1-\varphi\left[2^{-j+1}{|\tau|\over|\lambda|}\right]\right]=0$ if ${|\tau|\over |\lambda|}<2^{j-1}$ or ${|\tau|\over |\lambda|}>2^j$. On the other hand, 
$\p_\tau^\alpha\p_\lambda^\beta \varphi\left[2^j{|\tau|\over|\lambda|}\right]=0$ if ${|\tau|\over |\lambda|}<2^{-j}$ or ${|\tau|\over |\lambda|}>2^{-j+1}$.

A  direct computation shows
\bel{delta_j Diff Ineq}
\begin{array}{cc}\ds
\left| \p_\tau^\alpha\p_\lambda^\beta   \mathfrak{S}_j^+(\xi)\right|~\leq~\C_{\alpha~\beta}~ |\tau|^{-|\alpha|}\left[ 2^{-j}|\tau|\right]^{-|\beta|},
\\\\ \ds
\left| \p_\tau^\alpha\p_\lambda^\beta   \mathfrak{S}_j^-(\xi)\right|~\leq~\C_{\alpha~\beta}~ \left[2^{-j}|\lambda|\right]^{-|\alpha|}|\lambda|^{-|\beta|}
\end{array}
\eeq
for every multi-indices $\alpha,\beta$.

Define
\bel{F natural}
\begin{array}{cc}\ds
{^\natural}\F f(x)~=~\sum_{j>0} {^\natural}\F_j f(x)~=~\sum_{\ell\in\Z} {^\natural}\F_\ell f(x),
\\\\ \ds 
{^\natural}\F_j f(x)~=~\sum_{\ell\colon-j/2<\ell<j/2} \F_{\ell j} f(x),\qquad
{^\natural}\F_\ell f(x)~=~\sum_{j\colon j>2|\ell|} \F_{\ell j} f(x)
\end{array}
\eeq
together with
\bel{F sharp}
 \begin{array}{cc}\ds
  {^\sharp}\F f(x)~=~\sum_{j>0}{^\sharp}\F_j f(x),\qquad {^\sharp}\F_j f(x)~=~\int_{\R^\N} f(y)  {^\sharp}\Omega_j(x,y)dy, 
  \\\\ \ds
 {^\sharp}\Omega_j(x,y)~=~ \sum_{j/2\leq\ell<j} \Omega_{\ell j}(x,y)+{^+}\Omega_j(x,y),
 \\\\ \ds
{^+}\Omega_j(x,y)~=~\int_{\R^\N} e^{2\pi\i \big[\Phi(x,\xi)-y\cdot\xi \big]}\mathfrak{S}^+_j(\xi)\phi_j(\xi)\sigma(x,\xi)d\xi
\end{array}
\eeq
and
\bel{F flat}
 \begin{array}{cc}\ds
{^\flat}\F f(x)~=~\sum_{j>0}{^\flat}\F_j f(x),\qquad {^\flat}\F_j f(x)~=~\int_{\R^\N} f(y)  {^\flat}\Omega_j(x,y)dy, 
  \\\\ \ds
 {^\flat}\Omega_j(x,y)~=~ \sum_{-j<\ell\leq-j/2} \Omega_{\ell j}(x,y)+{^-}\Omega_j(x,y),
 \\\\ \ds
{^-}\Omega_j(x,y)~=~\int_{\R^\N} e^{2\pi\i \big[\Phi(x,\xi)-y\cdot\xi \big]}\mathfrak{S}^-_j(\xi)\phi_j(\xi)\sigma(x,\xi)d\xi.
\end{array}
\eeq
\begin{remark} Let $\F_o$ and ${^\natural}\F$,  ${^\sharp}\F$, ${^\flat}\F$ defined in (\ref{F_o}) and (\ref{F natural})-(\ref{F flat}) respectively. We find $\F=\F_o+{^\natural}\F+{^\sharp}\F+{^\flat}\F$.
\end{remark}
Denote
\bel{f mix-norm}
 \left\| f\right\|_{\L^{p_1}\L^{p_2}(\R^n\times\R^m)}~=~
\left\{\int_{\R^n}\left\{\int_{\R^m}\left|f(u,v)\right|^{p_2} dv\right\}^{p_1\over p_2} du\right\}^{1\over p_1},\qquad 1<p_1,~ p_2<\infty.
\eeq 

{\bf Proposition One}~~ {\it Suppose $\sigma\in\S^{-\rho}$ for $0\leq\rho\leq(\N-1)/2$ satisfying  (\ref{Class})-(\ref{rho_1, rho_2}). 
  We have
 
 {\bf (1)}
 \bel{L^2-result}
  \left\|\F_o f\right\|_{\L^2(\R^\N)}+\left\| {^\natural}\F f\right\|_{\L^2(\R^\N)}+ \left\| {^\sharp}\F f\right\|_{\L^2(\R^\N)}+ \left\| {^\flat}\F f\right\|_{\L^2(\R^\N)}~\leq~\C~\left\| f\right\|_{\L^2(\R^\N)},\qquad \sigma\in\S^0.
 \eeq 
 {\bf (2)}
\bel{F_t 2,p result}
\begin{array}{cc}\ds
 \left\|  {^\natural}\F_\ell f\right\|_{\L^2(\R^\N)}~\leq~\C_{\rho_1~\rho_2}~\left\| f\right\|_{\L^{p_1}\L^{p_2}(\R^n\times\R^m)},\qquad
 \ell\in\Z
\end{array}
\eeq
and
\bel{F sharp flat 2,p}
 \left\|  {^\sharp}\F f\right\|_{\L^2(\R^\N)}~\leq~\C_{\rho_1~\rho_2}~\left\| f\right\|_{\L^{p_1}\L^{p_2}(\R^n\times\R^m)},\qquad \left\|  {^\flat}\F f\right\|_{\L^2(\R^\N)}~\leq~\C_{\rho_1~\rho_2}~\left\| f\right\|_{\L^{p_1}\L^{p_2}(\R^n\times\R^m)} 
 \eeq
for
\bel{formula rho_1, rho_2}
{\rho_1\over n}~=~{1\over p_1}-{1\over 2},\qquad {\rho_2\over m}~=~{1\over p_2}-{1\over 2}.
\eeq
 {\bf (3)}
\bel{F*_t 2,p result}
\begin{array}{cc}\ds
 \left\| {^\natural} \F^*_\ell f\right\|_{\L^2(\R^\N)}~\leq~\C~2^{\big[{\N\over 2}-\rho-\ve\big]|\ell|} \left\| f\right\|_{\L^p(\R^\N)},\qquad
 {\rho\over \N}~=~{1\over p}-{1\over 2}
 \end{array}
\eeq
for some $\ve=\ve(\rho_1,\rho_2)>0$ and every $\ell\in\Z$. }\\
\begin{remark}
Given $0<\rho\leq(\N-1)/2$, we find 
$0<\rho_1<{n\over 2}$ and $0<\rho_2<{m\over 2}$
provided that $\rho_1, \rho_2$ belong to  (\ref{rho_1, rho_2}).
\end{remark}

\section{Proof of Theorem Two}
\setcounter{equation}{0}
Let $\F$ defined in (\ref{Ff})-(\ref{nondegeneracy}). We write
\bel{Ff Int}
\F f(x)~=~\int_{\R^\N} f(y) \left\{ \int_{\R^\N} e^{2\pi\i\big[\Phi(x,\xi)-y\cdot\xi\big]}\sigma(x,\xi)d\xi\right\} dy.
\eeq
Our main objective is to show
\bel{H^1 to L^1 EST}
\int_{\R^\N} \left| \F f(x)\right|dx~\leq~\C_{\rho_1~\rho_2}~\left\| f\right\|_{\H^1(\R^\N)},
\qquad\sigma\in\S^{-{\N-1\over2}}
\eeq
as well as  for $\F^*$. 
On the other hand, $\F$ is bounded on $\L^2(\R^\N)$ for $\sigma\in\S^0$ because  of the $\L^2$-estimate in (\ref{L^2-result}) and {\bf Remark 2.2}. Together, we  obtain the $\L^p$-norm inequality in (\ref{Result A})      by applying the complex interpolation theorem.  See {\bf 5.2}, chapter IV  of Stein \cite{Stein}.

Recall {\bf Remark 1.2}. $\H^1(\R^\N)$ has an atom decomposition as characterized by   Coifman \cite{Coifman}. 
Denote $a$ to be an $\H^1$-atom: $\supp a\subset\B_r(x_o)$, $|a(x)|\leq|\B_r(x_o)|^{-1}$ and $\ds\int_{\B_r(x_o)}a(x)dx=0$ where $\B(x_o)\subset\R^\N$ is a ball centered on some $x_o$ with radius $r>0$.  

By changing $y\mt y-x_o$ inside (\ref{Ff Int}), we find a new phase function $\Tilde{\Phi}(x,\xi)=\Phi(x,\xi)+x_o\cdot\xi$ which is again
real-valued, homogeneous of degree $1$ in $\xi$ and smooth in $(x,\xi)$ for $\xi\neq0$. Moreover, it satisfies   
 the non-degeneracy condition (\ref{nondegeneracy}).
Denote $\B_r\doteq\B_r(0)$. 
We assume $\supp a\subset \B_r$
  without lose of the generality.

To conclude (\ref{H^1 to L^1 EST}), it is suffice to prove 
\bel{H^1 est F_ta}
\begin{array}{lr}\ds
\int_{\R^\N} \left|\F a(x)\right|dx~\leq~\C,\qquad\sigma\in\S^{-{\N-1\over 2}}.
\end{array}
\eeq 
Let $\F_o$ and ${^\natural}\F, {^\sharp}\F, {^\flat}\F$  defined in (\ref{F_o}) and (\ref{F natural})-(\ref{F flat}). We find
\bel{F SUM}
\begin{array}{cc}\ds
\F a(x)~=~\F_o a(x)+{^\natural}\F a(x)+{^\sharp}\F a(x)+{^\flat}\F a(x),
\qquad 
{^\natural}\F a(x)~=~\sum_{\ell\in\Z}{^\natural}\F_\ell a(x).
\end{array}
\eeq
First,  $\F_o$ is $\L^1$-bounded as shown in (\ref{F_o L^1}). Next, 
each one of these partial operators has a compact support in $x$.  If $r\ge1$, we have
\bel{Est r>1}
\begin{array}{lr}\ds
\int_{\R^\N}\left|{^\natural}\F a(x)\right| dx
~\leq~\C~\left\| {^\natural}\F a\right\|_{\L^2(\R^\N)}\qquad\hbox{by Schwartz inequality}
\\\\ \ds~~~~~~~~~~~~~~~~~~~~~~~~~
~\leq~\C~\left\| a\right\|_{\L^2(\R^\N)}\qquad \hbox{\small{by (\ref{L^2-result})}}
\\\\ \ds~~~~~~~~~~~~~~~~~~~~~~~~~
~\leq~\C~r^{\N\big[-1+{1\over 2}\big]}\qquad\hbox{\small{( $|a(x)|\leq|\B_r|^{-1}$,~~$\supp a\subset\B_r$ )}}.
\end{array}
\eeq
Let $0<r<1$. 
Denote $\Q_{r\ell}\subset\R^\N$, so-called the {\it region of influence} associated to ${^\natural}\F_\ell$, satisfying
\bel{Q norm}
\left|\Q_{r\ell}\right|~\leq~\C~r\left\{\begin{array}{lr}\ds 2^{-m\ell},\qquad \ell\ge0,
\\ \ds 2^{n\ell},\qquad ~~\ell\leq0.
\end{array}
\right.
\eeq
This subset will be explicitly defined in the next section.

By using Schwartz inequality and (\ref{Q norm}), we find
\bel{Local est}
\begin{array}{lr}\ds
\int_{\Q_{r\ell}} \left|{^\natural}\F_\ell a(x)\right|dx~\leq~\C\left\|{^\natural}\F_\ell a\right\|_{\L^2(\R^\N)}~r^{1\over 2} \left\{\begin{array}{lr}\ds 2^{-m\ell/2},\qquad \ell\ge0,
\\ \ds 2^{n\ell/2},\qquad ~~\ell\leq0.
\end{array}\right.
\end{array}
\eeq
From (\ref{Local est}), by applying (\ref{F_t 2,p result}) in {\bf Proposition One} with $\rho={\N-1\over2}$,   we have
\bel{Local est 2}
\begin{array}{lr}\ds
\int_{\Q_{r\ell}} \left|{^\natural}\F_\ell a(x)\right|dx~\leq~\C_{\rho_1~\rho_2}\left\| a\right\|_{\L^{p_1}\L^{p_2}(\R^n\times\R^m)} 
r^{1\over 2}\left\{\begin{array}{lr}\ds 2^{-m\ell/2},\qquad \ell\ge0,
\\ \ds 2^{n\ell/2},\qquad ~~\ell\leq0
\end{array}\right.
\\\\ \ds~~~~~~~~~~~~~~~~~~~~~~~~~~
~\leq~\C_{\rho_1~\rho_2}~ r^{-\N}r^{\big[{n\over p_1}+{m\over p_2}\big]} r^{1\over 2} \left\{\begin{array}{lr}\ds 2^{-m\ell/2},\qquad \ell\ge0,
\\ \ds 2^{n\ell/2},\qquad ~~\ell\leq0
\end{array}\right.
\\\\ \ds~~~~~~~~~~~~~~~~~~~~~~~~~~
~=~\C_{\rho_1~\rho_2} \left\{\begin{array}{lr}\ds 2^{-m\ell/2},\qquad \ell\ge0,
\\ \ds 2^{n\ell/2},\qquad ~~\ell\leq0
\end{array}\right.
\end{array}
\eeq
where ${n\over p_1}+{m\over p_2}={n+m\over 2}+(\rho_1+\rho_2)={\N\over 2}+\rho={\N\over 2}+{\N-1\over 2}$.

In order to obtain (\ref{H^1 est F_ta}), we are left to show
\bel{Comple est l}
\begin{array}{lr}\ds
\int_{\R^\N\setminus\Q_{r\ell}} \left|{^\natural}\F_\ell a(x)\right|dx~\leq~\C\left\{\begin{array}{lr}\ds 2^{-m\ell/2},\qquad \ell\ge0,
\\ \ds 2^{n\ell/2},\qquad ~~\ell\leq0,
\end{array}\right.
\qquad 
\sigma\in\S^{-{\N-1\over 2}}
\end{array}
\eeq
and
\bel{Comple est}
\begin{array}{lr}\ds
\int_{\R^\N} \left|{^\sharp}\F a(x)\right|dx~\leq~\C_{\rho_1~\rho_2},\qquad
\int_{\R^\N} \left|{^\flat}\F a(x)\right|dx~\leq~\C_{\rho_1~\rho_2},
\qquad 
\sigma\in\S^{-{\N-1\over 2}}.
\end{array}
\eeq
\v
{\bf Proposition Two}~~{\it  Let $\Omega_{\ell j}$ and ${^\sharp}\Omega_j$, ${^\flat}\Omega_j$ for $\ell\in\Z, j>0$ defined in (\ref{Delta_lj F}) and (\ref{F sharp})-(\ref{F flat}) respectively.
Suppose $\sigma\in\S^{-{\N-1\over2}}$. We have
\bel{Result Size}
\begin{array}{cc}\ds
\int_{\R^\N} \left|\Omega_{\ell j}(x,y) \right|dx~\leq~\C \left\{\begin{array}{lr} 2^{-m\ell/2 },\qquad 0\leq\ell<j/2,
\\ \ds
2^{n\ell/2 },\qquad -j/2<\ell\leq0
\end{array}\right.
\end{array}
\eeq 
together with
\bel{Result Diff}
\begin{array}{ccc}\ds
\int_{\R^\N} \left|\Omega_{\ell j}(x,y) -\Omega_{\ell j}(x,0)\right|dx~\leq~\C~2^j|y| ~\left\{\begin{array}{lr} 2^{-m\ell/2},\qquad 0\leq\ell<j/2,
\\ \ds
2^{n\ell/2},\qquad -j/2<\ell\leq0
\end{array}\right.
\end{array}
\eeq 
and
\bel{Result Q}
\begin{array}{ccc}\ds
\int_{\R^\N\setminus \Q_{r\ell}} \left|\Omega_{\ell j}(x,y) \right|dx~\leq~\C~
{2^{-j}\over r} \left\{\begin{array}{lr} 2^{-m\ell/2},\qquad 0\leq\ell<j/2,
\\ \ds
2^{n\ell/2},\qquad -j/2<\ell\leq0
\end{array}\right.
\end{array}
\eeq 
for  $y\in\B_r$ whenever $2^j>r^{-1}$.
Moreover, we have
\bel{Result Size sharp+flat}
\int_{\R^\N} \left|{^\sharp}\Omega_j(x,y) \right|dx~\leq~\C ~2^{-\ve j},\qquad  \int_{\R^\N} \left|{^\flat}\Omega_j(x,y) \right|dx~\leq~\C~2^{-\ve j}
\eeq
for some $\ve=\ve(\rho_1,\rho_2)>0$.
}

Recall (\ref{Delta_lj F}) and (\ref{F natural})-(\ref{F flat}). We have
\bel{F natural l rewrite}
\begin{array}{cc}\ds
{^\natural}\F_\ell f(x)~=~\sum_{j\colon j>2|\ell|} \F_{\ell j}f(x),\qquad \F_{\ell j} f(x)~=~\int_{\R^\N} f(y) \Omega_{\ell j}(x,y)dy,
\\\\ \ds
{^\sharp}\F_\ell f(x)~=~\sum_{j>0} {^\sharp}\F_j f(x),\qquad {^\sharp}\F_j f(x)~=~\int_{\R^\N} f(y) {^\sharp}\Omega_j(x,y)dy,
\\\\ \ds
{^\flat}\F_\ell f(x)~=~\sum_{j>0} {^\sharp}\F_j f(x),\qquad {^\flat}\F_j f(x)~=~\int_{\R^\N} f(y) {^\flat}\Omega_j(x,y)dy.
\end{array}
\eeq
Let $a$ be an $\H^1$-atom where $\supp a\subset\B_r$ for $0<r<1$. For $2^j\leq r^{-1}$, we write
\bel{F_t cancella}
\int_{\R^\N} a(y)\Omega_{\ell j}(x,y)dy~=~\int_{\B_r} a(y) \left[\Omega_{\ell j}(x,y)-\Omega_{\ell j}(x,0)\right]dy
\eeq 
because 
$\ds\int_{\B_r} a(y)dy=0$. By using (\ref{Result Diff}) and (\ref{F_t cancella}), we find
\bel{Norm Est1}
\begin{array}{lr}\ds
\int_{\R^\N}\left|\int_{\R^\N} a(y)\Omega_{\ell j}(x,y)dy\right|dx
~\leq~\int_{\B_r} |a(y)|\left\{\int_{\R^\N} \left|\Omega_{\ell j}(x,y)-\Omega_{\ell j}(x,0)\right| dx\right\} dy
\\\\ \ds~~~~~~~~~~~~~~~~~~~~~~~~~~~~~~~~~~~~~~~~~~~~~~~
~\leq~ \C~ 2^j r\left\{\begin{array}{lr} 2^{-m\ell/2},\qquad 0\leq\ell<j/2,
\\ \ds
2^{n\ell/2},\qquad -j/2<\ell<0.
\end{array}\right.
\end{array}
\eeq
By summing over every $j>2|\ell|$ and $2^j\leq r^{-1}$, we find
\bel{Sum1}
\begin{array}{lr}\ds
 \sum_{j\colon j>2|\ell|,~2^j\leq r^{-1}} \int_{\R^\N}\left|\int_{\R^\N} a(y)\Omega_{\ell j}(x,y)dy\right|dx
~\leq~\C \left\{\begin{array}{lr} 2^{-m\ell/2},\qquad \ell\ge0,
\\ \ds
2^{n\ell/2},\qquad~~~ \ell\leq0.
\end{array}\right.
\end{array}
\eeq
On the other hand, for $2^j>r^{-1}$,  we have
\bel{Norm Est2}
\begin{array}{lr}\ds
\int_{\R^\N\setminus\Q_{r\ell}}\left|\int_{\R^\N} a(y)\Omega_{\ell j}(x,y)dy\right|dx
~\leq~\int_{\B_r} |a(y)|\left\{\int_{\R^\N\setminus\Q_{r\ell}} \left|\Omega_{\ell j}(x,y)\right| dx\right\} dy
\\\\ \ds~~~~~~~~~~~~~~~~~~~~~~~~~~~~~~~~~~~~~~~~~~~~~~~~~~~~
~\leq~\C~{2^{-j}\over r}\left\{\begin{array}{lr} 2^{-m\ell/2},\qquad 0\leq\ell<j/2,
\\ \ds
2^{n\ell/2},\qquad -j/2<\ell\leq0
\end{array}\right.\qquad \hbox{\small{by (\ref{Result Q})}}.
\end{array}
\eeq
By summing over every  $j>2|\ell|$  and $2^j>r^{-1}$, we find
\bel{Sum2}
\begin{array}{lr}\ds
\sum_{j\colon j>2|\ell|,~2^j>r^{-1}} \int_{\R^\N\setminus\Q_{r\ell}}\left|\int_{\R^\N} a(y)\Omega_{\ell j}(x,y)dy\right|dx
~\leq~\C~\left\{\begin{array}{lr} 2^{-m\ell/2},\qquad \ell\ge0,
\\ \ds
2^{n\ell/2},\qquad ~~~\ell\leq0.
\end{array}\right. 
\end{array}
\eeq
By putting together (\ref{Sum1}) and (\ref{Sum2}), we conclude (\ref{Comple est l}).

Furthermore, by applying (\ref{Result Size sharp+flat}), we have
\bel{INT EST sharp+flat}
\begin{array}{lr}\ds
\left\| \F^\sharp_j a\right\|_{\L^1(\R^\N)}~\leq~\int_{\R^\N} |a(y)| \left\{\int_{\R^\N} \left|{^\sharp}\Omega_j(x,y)\right|dx\right\} dy
~\leq~\C~2^{-\ve j},
\\\\ \ds
\left\| \F^\flat_j a\right\|_{\L^1(\R^\N)}~\leq~\int_{\R^\N} |a(y)| \left\{\int_{\R^\N} \left|{^\flat}\Omega_j(x,y)\right|dx\right\} dy
~\leq~\C~2^{-\ve j}
\end{array}
\eeq
for some $\ve=\ve(\rho_1,\rho_2)>0$. 
Clearly, (\ref{INT EST sharp+flat}) implies (\ref{Comple est}) by summing all $j>0$.

\section{A second  dyadic decomposition}
\setcounter{equation}{0}
Given $j>0$ and $\ell\in\Z$, 
our assertion is split into two, $w.r.t$
\bel{Cases}
\hbox{\bf Case One}: ~~-j/2<\ell<j/2,\qquad  \hbox{\bf Case Two}: ~~\ell\leq -j/2~~~\hbox{and}~~~\ell\ge j/2.
\eeq

{\bf Case One}~~ Recall (\ref{Delta_lj F}) and (\ref{F natural}). The kernel of $\F_{\ell j}$  is  given by
\bel{Omega_lj} 
\begin{array}{cc}\ds
\Omega_{\ell j}(x,y)~=~
\int_{\R^\N} e^{2\pi\i \big[\Phi(x,\xi)-y\cdot\xi \big]}\delta_\ell(\xi)\phi_j(\xi)\sigma(x,\xi)d\xi
\end{array}
\eeq
where $\delta_\ell$ and $\phi_j$ are defined in (\ref{delta_t}) and (\ref{phi_j}). 

Denote $\left\{ \xi^\nu_j\right\}_\nu$ to be a collection of points  equally distributed on  $\mathds{S}^{\N-1}$ with a grid length between $2^{-j/2-1}$ and $2^{-j/2}$. 
Let
\bel{Gamma_j}
\Gamma_j^\nu~=~\Bigg\{\xi\in\R^\N~\colon \left|{\xi\over |\xi|}-\xi^\nu_j\right|\leq 2^{-j/2+1}\Bigg\}.
\eeq
Let $\varphi\in\mathcal{C}^\infty_o(\R)$ such that $\varphi(t)=1$ if $|t|\leq1$ and $\varphi(t)=0$ for $|t|>2$.
  We define
\bel{phi^v_j}
\varphi^\nu_j(\xi)~=~{\ds\varphi\Bigg[2^{j/2}\left|{\xi\over |\xi|}-\xi^{\nu}_{j}\right|\Bigg]\over \ds\sum_{\nu\colon\{\xi^\nu_j\}\subset\mathds{S}^{\N-1}} \varphi\Bigg[2^{j/2}\left|{\xi\over |\xi|}-\xi^{\nu}_{j}\right|\Bigg]}.
\eeq
Observe that $\supp\varphi^\nu_j\subset\Gamma^\nu_j$. On the other hand, $\supp\delta_\ell\subset\Lambda_\ell$ defined in (\ref{Cone}). Therefore, it is essential to assert
$\Gamma^\nu_j\cap\Lambda_\ell\neq\emptyset$.

From (\ref{Omega_lj}) and (\ref{phi^v_j}), we write
\bel{Omega_lj sum}
\begin{array}{cc}\ds
\Omega_{\ell j}(x,y)~=~\sum_{\nu\colon\{\xi^\nu_j\}\subset\mathds{S}^{\N-1}} ~\Omega_{\ell j}^\nu(x,y),
\\\\ \ds
\Omega_{\ell j}^\nu(x,y)~=~
\int_{\R^\N} e^{2\pi\i \big[\Phi(x,\xi)-y\cdot \xi \big]}\delta_\ell(\xi)\phi_j(\xi)\varphi^\nu_j(\xi)\sigma(x,\xi)d\xi.
\end{array}
\eeq
\begin{remark} As a geometric fact, there are at most a constant multiple of 
\bel{terms number}
2^{j\big({\N-1\over 2}\big)}\left\{\begin{array}{lr}\ds 2^{-\ell m},\qquad \ell\ge0,
\\ \ds
2^{n\ell},\qquad~~~ \ell\leq0
\end{array}
\right.
\eeq
many $\nu$'s such that $\xi^\nu_j\in\mathds{S}^{\N-1}\cap\Lambda_\ell$.
\end{remark}
Let $\c>0$  be some large constant depending on $\Phi$. Given $\xi^\nu_j\in\mathds{S}^{\N-1}$, we consider
\bel{R xi natural}
\hbox{\bf R}^\nu_j~=~\left\{x\in\supp\sigma~\colon \left|\left(\nabla_\xi\Phi\right)\left(x,\xi^\nu_j\right)\cdot \xi^\nu_j\right|\leq\c 2^{-j},~\left|\left(\nabla_\xi \Phi\right)\left(x,\xi^\nu_j\right)\right|\leq\c 2^{-j/2}\right\}.
\eeq
The region of influence  $\Q_{r\ell}\subset\R^\N$, $0<r<1$ associated to ${^\natural}\F_\ell$ is defined by
  \bel{Q natural}
\Q_{r\ell}~=~\bigcup_{j\colon2^{-j}\leq r} ~\bigcup_{\nu\colon\xi_j^\nu\in\mathds{S}^{\N-1}\cap\Lambda_\ell} \hbox{\bf R}^\nu_j.
\eeq
Each $\hbox{\bf R}^\nu_j$ has a volume bounded by $\C \c^\N2^{-j}2^{-j\left({\N-1\over 2}\right)}$. Recall {\bf Remark 4.1}. A direct computation shows 
\bel{Q natural norm}
\left|\Q_{r\ell}\right|~\leq~\C~r\left\{\begin{array}{lr}\ds 2^{-\ell m},\qquad \ell\ge0,
\\ \ds
2^{n\ell},\qquad~~~ \ell\leq0.
\end{array}
\right.
\eeq

{\bf Case Two}~~Recall (\ref{F sharp}) and (\ref{F flat}). The kernel of ${^\sharp}\F_j$ is given by
\bel{Omega_j+} 
\begin{array}{lr}\ds
{^\sharp}\Omega_j(x,y)~=~\sum_{j/2\leq\ell<j} \Omega_{\ell j}(x,y)+{^+}\Omega_j(x,y)
\\\\ \ds~~~~~~~~~~~~~~
~=~ \sum_{j/2\leq\ell<j}\int_{\R^\N} e^{2\pi\i \big[\Phi(x,\xi)-y\cdot\xi \big]}\delta_\ell(\xi)\phi_j(\xi)\sigma(x,\xi)d\xi
\\\\ \ds~~~~~~~~~~~~~~
~+~ \int_{\R^\N} e^{2\pi\i \big[\Phi(x,\xi)-y\cdot\xi \big]}\mathfrak{S}_j^+(\xi)\phi_j(\xi)\sigma(x,\xi)d\xi
\end{array}
\eeq
where $\mathfrak{S}_j^+$ is defined in (\ref{delta_j+-}).

Write $\xi=(\tau,\lambda)\in\R^n\times\R^m$.  $\left\{ \tau^\nu_j\right\}_\nu$ is a collection of points  equally distributed on  $\mathds{S}^{n-1}$ with a grid length between $2^{-j/2-1}$ and $2^{-j/2}$.

Denote $\xi^\nu_j=(\tau^\nu_j,0)\in\R^n\times\R^m$. We consider
\bel{phi^v_j +}
{^+}\varphi^\nu_j(\xi)~=~{\ds\varphi\Bigg[2^{j/2}\left|{\xi\over |\xi|}-\xi^{\nu}_{j}\right|\Bigg]\over \ds\sum_{\nu\colon\{\tau^\nu_j\}\subset\mathds{S}^{n-1}} \varphi\Bigg[2^{j/2}\left|{\xi\over |\xi|}-\xi^{\nu}_{j}\right|\Bigg]},\qquad \supp {^+}\varphi^\nu_j\subset\Gamma^\nu_j.
\eeq
From (\ref{Omega_j+}), we write
\bel{Omega_j sum+}
\begin{array}{ccc}\ds
\Omega_{\ell j}(x,y)~=~\sum_{\nu\colon\{ \tau^\nu_j\}_\nu\subset\mathds{S}^{n-1}} {^+}\Omega^\nu_{\ell j}(x,y),
\\\\ \ds
{^+}\Omega^\nu_{\ell j}(x,y)~=~\int_{\R^\N} e^{2\pi\i \big[\Phi(x,\xi)-y\cdot\xi \big]}\delta_\ell(\xi)\phi_j(\xi){^+}\varphi^\nu_j(\xi)\sigma(x,\xi)d\xi;
\\\\ \ds
{^+}\Omega_j(x,y)~=~\sum_{\nu\colon\{ \tau^\nu_j\}_\nu\subset\mathds{S}^{n-1}} {^+}\Omega^\nu_j(x,y),
\\\\ \ds
{^+}\Omega^\nu_j(x,y)~=~
\int_{\R^\N} e^{2\pi\i \big[\Phi(x,\xi)-y\cdot\xi \big]}\mathfrak{S}_j^+(\xi)\phi_j(\xi){^+}\varphi^\nu_j(\xi)\sigma(x,\xi)d\xi.
\end{array}
\eeq

\begin{remark}
$\ds\Cup_{j/2\leq\ell<j} \supp\delta_\ell\cup
 \supp\mathfrak{S}^+_j \subset \bigcup_{\nu\colon~\{ \tau^\nu_j\}_\nu\subset\mathds{S}^{n-1}} \Gamma^\nu_j$.
\end{remark}
\begin{remark} There are at most a constant multiple of $2^{j\big({n-1\over 2}\big)}$ many $\nu$'s in $\left\{ \tau^\nu_j\right\}_\nu$.
\end{remark}
On the other hand, the kernel of ${^\flat}\F_j$ is given by
\bel{Omega_j-} 
\begin{array}{lr}\ds
{^\flat}\Omega_j(x,y)~=~\sum_{-j<\ell\leq-j/2} \Omega_{\ell j}(x,y)+{^-}\Omega_j(x,y)
\\\\ \ds~~~~~~~~~~~~~~
~=~ \sum_{-j<\ell\leq-j/2}\int_{\R^\N} e^{2\pi\i \big[\Phi(x,\xi)-y\cdot\xi \big]}\delta_\ell(\xi)\phi_j(\xi)\sigma(x,\xi)d\xi
\\\\ \ds~~~~~~~~~~~~~~
~+~ \int_{\R^\N} e^{2\pi\i \big[\Phi(x,\xi)-y\cdot\xi \big]}\mathfrak{S}_j^-(\xi)\phi_j(\xi)\sigma(x,\xi)d\xi
\end{array}
\eeq
where $\mathfrak{S}_j^-$ is defined in (\ref{delta_j+-}).

Write $\xi=(\tau,\lambda)\in\R^n\times\R^m$.  $\left\{ \lambda^\nu_j\right\}_\nu$ is a collection of points  equally distributed on  $\mathds{S}^{m-1}$ with a grid length between $2^{-j/2-1}$ and $2^{-j/2}$.
Denote $\xi^\nu_j=(0,\lambda^\nu_j)\in\R^n\times\R^m$. We consider
\bel{phi^v_j -}
{^-}\varphi^\nu_j(\xi)~=~{\ds\varphi\Bigg[2^{j/2}\left|{\xi\over |\xi|}-\xi^{\nu}_{j}\right|\Bigg]\over \ds\sum_{\nu\colon\{\lambda^\nu_j\}\subset\mathds{S}^{m-1}} \varphi\Bigg[2^{j/2}\left|{\xi\over |\xi|}-\xi^{\nu}_{j}\right|\Bigg]},\qquad \supp {^-}\varphi^\nu_j\subset\Gamma^\nu_j.
\eeq
From (\ref{Omega_j-}), we write
\bel{Omega_j sum-}
\begin{array}{ccc}\ds
\Omega_{\ell j}(x,y)~=~\sum_{\nu\colon\{ \lambda^\nu_j\}_\nu\subset\mathds{S}^{m-1}} {^-}\Omega^\nu_{\ell j}(x,y),
\\\\ \ds
{^-}\Omega^\nu_{\ell j}(x,y)~=~\int_{\R^\N} e^{2\pi\i \big[\Phi(x,\xi)-y\cdot\xi \big]}\delta_\ell(\xi)\phi_j(\xi){^-}\varphi^\nu_j(\xi)\sigma(x,\xi)d\xi;
\\\\ \ds
{^-}\Omega_j(x,y)~=~\sum_{\nu\colon\{ \lambda^\nu_j\}_\nu\subset\mathds{S}^{m-1}} {^-}\Omega^\nu_j(x,y),
\\\\ \ds
{^-}\Omega^\nu_j(x,y)~=~
\int_{\R^\N} e^{2\pi\i \big[\Phi(x,\xi)-y\cdot\xi \big]}\mathfrak{S}_j^-(\xi)\phi_j(\xi){^+}\varphi^\nu_j(\xi)\sigma(x,\xi)d\xi.
\end{array}
\eeq
\begin{remark}
$\ds\Cup_{-j<\ell\leq-j/2}\supp\delta_\ell\cup
 \supp\mathfrak{S}^-_j \subset \bigcup_{\nu\colon~\{ \lambda^\nu_j\}_\nu\subset\mathds{S}^{m-1}} \Gamma^\nu_j$.
\end{remark}
\begin{remark} There are at most a constant multiple of $2^{j\big({m-1\over 2}\big)}$ many $\nu$'s in $\left\{ \lambda^\nu_j\right\}_\nu$.
\end{remark}

\section{Proof of Proposition Two}
\setcounter{equation}{0}
By symmetry, we prove the regarding estimates in (\ref{Result Size})-(\ref{Result Size sharp+flat}) for  $\ell\ge0$ only.

{\bf Case One}~~Let $0\leq\ell<j/2$. Recall $\left\{ \xi^\nu_j\right\}_\nu$ which is a collection of points  equally distributed on  $\mathds{S}^{\N-1}$ with a grid length between $2^{-j/2-1}$ and $2^{-j/2}$. 

Given $\nu$ fixed, we assert
$\xi=\L_\nu \eta$ where $\L_\nu$ is an $\N\times \N$-orthogonal matrix   with $\det\L_\nu=1$. 
Write $\eta=(\eta_1,\eta')\in\R\times\R^{\N-1}$.  Denote $\eta^1=(\eta_1,0)\in\R\times\R^{\N-1}$. We require that $\eta^1$ is in the same direction of $\xi^\nu_j$. 

From (\ref{Omega_lj sum}), we have
\bel{Omega_lj v}
\begin{array}{cc}\ds
\Omega_{\ell j}(x,y)~=~\sum_{\nu\colon\{\xi^\nu_j\}\subset\mathds{S}^{\N-1}} ~\Omega_{\ell j}^\nu(x,y),
\\\\ \ds
\Omega_{\ell j}^\nu(x,y)~=~
\int_{\R^\N} e^{2\pi\i \big[\Phi(x,\L_\nu\eta)-y\cdot \L_\nu \eta \big]}\delta_\ell(\L_\nu\eta)\phi_j(\L_\nu\eta)\varphi^\nu_j(\L_\nu\eta)\sigma(x,\L_\nu\eta)d\eta
\end{array}
\eeq
where $\delta_\ell$, $\phi_j$ and $\varphi_j^\nu$ are  defined in (\ref{delta_t}), (\ref{phi_j}) and (\ref{phi^v_j}) respectively.

Denote $\eta^\nu_j=(1,0)\in\R\times\R^{\N-1}$ so that $\xi^\nu_j=\L_\nu \eta^\nu_j$. 
Consider
\bel{Phi split}
\begin{array}{cc}\ds
\Phi_\nu(x,\eta)~=~\Phi(x,\L_\nu \eta),
\\\\ \ds
\Phi_\nu(x,\eta)-y\cdot \L_\nu\eta
~=~\left[\Big(\nabla_\eta\Phi_\nu\Big)\left(x,\eta_j^\nu\right)-\L_\nu^T y\right]\cdot\eta~+~\Psi_\nu(x,\eta),
\\\\ \ds
\Psi_\nu(x,\eta)~=~\Phi_\nu(x,\eta)-\Big(\nabla_\eta\Phi_\nu\Big)\left(x,\eta_j^\nu\right)\cdot\eta.
\end{array}
\eeq
From (\ref{Omega_lj v})-(\ref{Phi split}), we rewrite
\bel{Theta}\begin{array}{cc}\ds
\Omega^\nu_{\ell j}(x,y)~=~\int_{\R^n} e^{2\pi\i\Big[\big(\nabla_\eta\Phi_\nu\big)\big(x,\eta_j^\nu\big)-\L_\nu^T y\Big]\cdot\eta}\Theta^\nu_{\ell j}(x,\eta)d\eta,
\\\\ \ds
\Theta^\nu_{\ell j}(x,\eta)~=~e^{2\pi\i\Psi_\nu(x,\eta)}\delta_\ell(\L_\nu\eta)\phi_j(\L_\nu\eta)\varphi^\nu_j(\L_\nu\eta)\sigma(x,\L_\nu\eta).
\end{array}
\eeq
Define a differential operator
\bel{D natural}
\D_j~=~I+2^{2j} \partial_{\eta_1}^2+2^j  \Delta_{\eta'},\qquad j>0.
\eeq
\begin{lemma} Let $\sigma\in\S^{-{\N-1\over 2}}$.
Given $j>0$ and $0\leq\ell<j/2$, we have
\bel{Theta Diff Ineq}
\left|\D_j^N \Theta^\nu_{\ell j}(x,\eta)\right|~\leq~\C_N~2^{-j\left({\N-1\over 2}\right)}2^{\ell\big({m\over 2}\big)},\qquad N\ge0.
\eeq
\end{lemma}
{\bf Proof}~~Recall $\supp\varphi^\nu_j\subset\Gamma^\nu_j$ defined in (\ref{Gamma_j}). 
We find
 \bel{eta size}
2^{j-1}\leq|\eta_1|\leq2^{j+1},\qquad |\eta'|\leq \C 2^{j/2}
\eeq
whenever $\L_\nu \eta\in\Gamma_j^\nu\cap \left\{2^{j-1}\leq|\eta|\leq2^{j+1}\right\}$.

We claim
\bel{d Est Psi}
\left|\p_\eta^\alpha\Psi(x,\eta)\right|~\leq~\C_\alpha~2^{-j\alpha_1}2^{-(j/2)\big[|\alpha|-\alpha_1\big]}
\eeq
for every multi-index $\alpha$. 

To show (\ref{d Est Psi}),  we follow the lines in
 {\bf 4.5}, chapter IX of Stein \cite{Stein}.

Observe that $\Psi_\nu(x,\eta)$ in (\ref{Phi split}) is homogeneous of degree $1$ in $\eta$: $\Psi_\nu(x,r\eta)=r\Psi_\nu(x,\eta),r>0$. For every multi-index $\alpha$, we have
\bel{Psi est1}
\left|\partial_\eta^\alpha \Psi_\nu(x,\eta)\right|~\leq~\C_\alpha~|\eta|^{1-|\alpha|}.
\eeq
Indeed, denote $\vartheta=r\eta$. We have $\partial_\eta^\alpha \Psi_\nu(x,r\eta)=r^{|\alpha|}\partial_\vartheta^\alpha\Psi_\nu(x,\vartheta)=r\partial_\eta^\alpha\Psi_\nu(x,\eta)$. Choose $r=|\eta|^{-1}$. We find $\partial_\eta^\alpha\Psi_\nu(x,\eta)=|\eta|^{1-|\alpha|} \partial_\vartheta^\alpha\Psi_\nu(x,\vartheta)$ where $|\vartheta|=1$.

Suppose $|\alpha|-\alpha_1\ge2$. From (\ref{eta size}) and (\ref{Psi est1}), we have
\bel{Psi est2}
\begin{array}{lr}\ds
\left|\partial_\eta^\alpha \Psi_\nu(x,\eta)\right|~\leq~\C_\alpha~2^{-j\big(|\alpha|-1\big)}~\leq~\C_\alpha~2^{-j\alpha_1}2^{-j\big[|\alpha|-\alpha_1-1\big]} 
\\\\ \ds~~~~~~~~~~~~~~~~~~~~~~~~~~~~~~~~~~~~~~~~~~~~~~~
~\leq~\C_\alpha~2^{-j\alpha_1}2^{-(j/2)\big[|\alpha|-\alpha_1\big]}. 
\end{array}
\eeq
Consider $|\alpha|-\alpha_1=0$ or $1$. Recall  that $\Phi(x,\xi)$ is homogeneous of degree $1$ in $\xi$. Write $\ds\Phi(x,\xi)=\int_0^1 {d\over dr}\left[\Phi(x,r\xi)\right] dr$. By integration by parts $w.r.t~ r$ and taking into account that $\Phi(x,r\xi)=r\Phi(x,\xi)$, we find
$\Phi(x,\xi)=\nabla_\xi\Phi(x,\xi)\cdot\xi$. This together with the identity 
 \bel{xi eta iden}
 \L_\nu^T \nabla_\xi \Phi(x,\xi)~=~\nabla_\eta \left[\Phi(x,\L_\nu \eta)\right],\qquad \xi=\L_\nu \eta
 \eeq
imply
 \bel{Psi rewrite}
 \Psi_\nu(x,\eta)~=~\nabla_\eta\Phi_\nu(x,\eta)\cdot\eta-\Big(\nabla_\eta\Phi_\nu\Big)\left(x,\eta_j^\nu\right)\cdot\eta.
 \eeq 
 Note that $(\nabla_\eta\Phi_\nu)(x,\eta^1)=(\nabla_\eta\Phi_\nu)(x,\eta^\nu_j)$ because $\nabla_\eta\Phi_\nu(x, \eta)$ is homogeneous of degree $0$ in $\eta$.  
 Therefore, we find $\partial_{\eta_1}\Psi_\nu(x,\eta^1)=0$.
 A direct computation shows
 \bel{p_eta_1 N}
 \begin{array}{lr}\ds
 ~~~~~~~~~~~~~~~~~~~\left( \nabla_\eta \Psi_\nu\right)(x,\eta^1)~=~0;
 \\\\ \ds 
 0~=~\Big(\partial_{\eta_1}^N \nabla_{\eta'}\Psi_\nu\Big)(x,\eta^1)
 ~=~ \Big(\nabla_{\eta'}\partial_{\eta_1}^N \Psi_\nu\Big)(x,\eta^1),\qquad N\ge0.
 \end{array}
 \eeq 
By writing out  the Taylor expansion of $\partial_{\eta_1}^N \Psi_\nu(x,\eta)$ and $\nabla_{\eta'}\partial_{\eta_1}^N \Psi_\nu(x,\eta)$ respectively in the $\eta'$-subspace and using (\ref{p_eta_1 N}), we find
 \bel{p_eta_1 N norm}
 \partial_{\eta_1}^N \Psi_\nu(x,\eta)~=~ \O\left(|\eta'|^2|\eta|^{-N-1}\right),\qquad  
 \nabla_{\eta'} \partial_{\eta_1}^N \Psi_\nu(x,\eta)~=~ \O\left(|\eta'||\eta|^{-N-1}\right)
  \eeq
    simultaneously for every $N\ge0$.
  
 Recall (\ref{eta size}). We further have 
 \bel{p_eta_1 N Est}
  \left|\partial_{\eta_1}^N \Psi_\nu(x,\eta)\right|~\leq~\C_N  2^{-j N},\qquad \left|\nabla_{\eta'}\partial_{\eta_1}^N \Psi_\nu(x,\eta)\right|~\leq~\C_N  2^{-j N}2^{-j/2},\qquad N\ge0.
  \eeq
By putting together (\ref{Psi est2}) and (\ref{p_eta_1 N Est}),  we obtain (\ref{d Est Psi}). 

Next, we prove the differential inequality in (\ref{Theta Diff Ineq}) by showing
\bel{Ineq l}
\begin{array}{lr}
\left|\p_\eta^\alpha \left[\delta_\ell (\L_\nu \eta)\phi_j(\L_\nu\eta)\varphi_j^\nu \left(\L_\nu\eta\right)\sigma(x,\L_\nu \eta)\right] \right|
\\\\ \ds
~\leq~\C_\alpha~2^{-j\left({\N-1\over 2}\right)}2^{\ell\big({m\over 2}\big)} 2^{-j \alpha_1} 2^{-(j/2)\big[|\alpha|-\alpha_1\big]}
\end{array}
\eeq
for every multi-index $\alpha$.

Recall $\varphi^\nu_j$  defined in  (\ref{phi^v_j}) where $\supp\varphi_j^\nu\subset\Gamma^\nu_j$ is given in (\ref{Gamma_j}).
Geometrically, for $\L_\nu\eta\in\Gamma^\nu_j$, the angle between $\eta$ and $\eta^1$ is bounded by $ \arcsin\Big(2\cdot2^{-j/2}\Big)$.

Denote $\rho=|\eta|$. By using polar coordinates, we find
\bel{radial derivative}
{\p\over \p\eta_1}~=~\left({\p \rho\over \p \eta_1}\right){\p\over \p \rho}+\O\left(2^{-j/2}\right)\cdot\nabla_{\eta'}.
\eeq
Because $\varphi^\nu_j\left(\L_\nu \eta\right)$  is homogeneous of degree zero in $\eta$,
 we find
\bel{partial 0}
\p_\rho\left[\varphi^\nu_j(\L_\nu \eta)\right]~=~0.
\eeq
A direct computation shows
\bel{nu Diff est}
\begin{array}{cc}\ds
\left|\p_\eta^\alpha \left[\varphi_j^{\nu}\left(\L_\nu\eta\right)\right]\right|~\leq~\C_{\alpha} ~2^{|\alphaup| j/2}|\eta|^{-|\alphaup|}
\end{array}
\eeq
for every multi-index $\alpha$. 
From  (\ref{radial derivative})-(\ref{nu Diff est}), we conclude
\bel{nu Diff est 1}
\begin{array}{cc}\ds
\left|\p_{\eta_1}^N \left[\varphi_j^\nu \left(\L_\nu\eta\right)\right]\right|~\leq~\C_{N} ~|\eta|^{-N},\qquad N\ge0.
\end{array}
\eeq
Recall that $\delta_\ell$ satisfies the differential inequality in (\ref{delta_l Diff Ineq}).
Moreover, it has a support inside the dyadic cone $\Lambda_\ell$ defined in (\ref{Cone}). 
Let $\xi=(\tau,\lambda)\in\R^n\times\R^m$.    
We find $|\tau|\sim2^j$,  $|\lambda|\sim2^{j-\ell}$ if $\xi\in \Lambda_\ell\cap\left\{2^{j-1}\leq|\xi|\leq2^{j+1}\right\}$.

Let $\sigma\in\S^{-{\N-1\over2}}$ satisfying the differential inequality in (\ref{Class})-(\ref{rho_1, rho_2}). Note that  $|\tau|\sim2^\ell|\lambda|$, $\left({1\over1+|\tau|}\right)^{\rho_1}\left({1\over 1+|\lambda|}\right)^{\rho_2}\lesssim\left({1\over1+|\tau|}\right)^{n-1\over 2}\left({1\over 1+|\lambda|}\right)^{m\over 2}$  for which $\rho_1> {n-1\over 2}$, $\rho_2>{m-1\over 2}$ and  $\rho_1+\rho_2={\N-1\over 2}$.
We have
\bel{INEQ l}
\begin{array}{lr}\ds
\left| \partial^\alpha_\tau\partial^\beta_\lambda \left[\delta_\ell(\xi) \phi_j(\xi) \sigma(x,\xi)\right]\right|
~\leq~\C_{\alpha~\beta} \left({1\over1+|\tau|}\right)^{n-1\over 2}\left({1\over 1+|\lambda|}\right)^{m\over 2} |\tau|^{-|\alpha|}|\lambda|^{-|\beta|}
\\\\ \ds~~~~~~~~~~~~~~~~~~~~~~~~~~~~~~~~~~~~~~~~~~
~\leq~\C_{\alpha~\beta}~2^{-j\big({\N-1\over 2}\big)}2^{\ell\big({m\over 2}\big)} 2^{-j |\alpha|} 2^{-(j-\ell)|\beta|}
\end{array}
\eeq
for every multi-indices $\alpha,\beta$.

Consider $\xi=\L_\nu\eta$. Write out
\bel{tau_i,lambda_j}
\begin{array}{cc}\ds
\tau_i~=~a_{i1}\eta_1+\O(1)\cdot \eta',\qquad i=1,\ldots,n,
\\\\ \ds
\lambda_k~=~a_{k1}\eta_1+\O(1)\cdot \eta',\qquad k=n+1,\ldots,n+m
\end{array}
\eeq
where $a_{i1}, i=1,\ldots,n$ is the entry on the $i$-th row and the first column of $\L_\nu$.  Vice versa for $a_{k1}, k=n+1,\ldots,n+m$. 

Suppose $\L_\nu\eta\in\Gamma^\nu_j\cap\Lambda_\ell\cap\{2^{j-1}\leq|\eta|\leq2^{j+1}\}$.
By putting together (\ref{eta size}) and (\ref{tau_i,lambda_j}), we find
\bel{a_k1}
|a_{k1}|~\leq~\C2^{-\ell},\qquad k~=~n+1,\ldots,n+m.
\eeq
From (\ref{INEQ l}) and (\ref{tau_i,lambda_j})-(\ref{a_k1}), we have
\bel{INEQ eta}
\begin{array}{lr}\ds
\left| \partial^\alpha_\eta \left[\delta_\ell(\L_\nu\eta) \phi_j(\L_\nu\eta) \sigma(x,\L_\nu\eta)\right]\right|
~\leq~\C_{\alpha~\beta}~2^{-j\big({\N-1\over 2}\big)}2^{\ell\big({m\over 2}\big)} 2^{-j \alpha_1} 2^{-(j/2)\big[|\alpha|-\alpha_1\big]}
\end{array}
\eeq
for every multi-index $\alpha$. 

By combining (\ref{INEQ eta})
 with (\ref{nu Diff est})-(\ref{nu Diff est 1}), we conclude (\ref{Ineq l}) as desired. \endproof

We finish the following estimates within 3 steps. 

{\bf 1.} ~~Let $\Omega^\nu_{\ell j}(x,y)$ and $\Theta^\nu_{\ell j}(x,\eta)$ defined in (\ref{Theta}).
In particular, we find 
\bel{Theta supp}
\left|\supp \Theta^\nu_{\ell j}\right|\leq\C2^j 2^{j\left({\N-1\over 2}\right)}.
\eeq
By applying (\ref{Theta Diff Ineq}), an $N$-fold integration by parts associated to $\D_j$ defined in (\ref{D natural}) shows 
\bel{Omega lj v size}
\begin{array}{lr}\ds
\left|\Omega^\nu_{\ell j}(x,y)\right|~\leq~
\\ \ds
\left\{1+4\pi^2 2^{2j}\left[\Big(\nabla_\eta\Phi_\nu\Big)\left(x,\eta_j^\nu\right)-\L_\nu^T y\right]_1^2+4\pi^2 2^{j}\sum_{i=2}^\N \left[\Big(\nabla_\eta\Phi_\nu\Big)\left(x,\eta_j^\nu\right)-\L_\nu^T y\right]_i^2 \right\}^{-N}\int_{\R^\N} \left|\D_j^N \Theta^\nu_{\ell j}(x,\eta)\right| d\eta
\\\\ \ds
~\leq~\C_{N}~2^{-j\left({\N-1\over 2}\right)}2^{\ell\big({m\over 2}\big)}~2^j 2^{j\left({\N-1\over 2}\right)}
\\ \ds~~~~~~~
\left\{1+4\pi^2 2^{2j}\left[\Big(\nabla_\eta\Phi_\nu\Big)\left(x,\eta_j^\nu\right)-\L_\nu^T y\right]_1^2+4\pi^2 2^{j}\sum_{i=2}^\N \left[\Big(\nabla_\eta\Phi_\nu\Big)\left(x,\eta_j^\nu\right)-\L_\nu^T y\right]_i^2 \right\}^{-N}.
\end{array}
\eeq
Consider a local diffeomorphism  
\bel{local diffeo}
\mathcal{X}_\Phi~\colon~ x~\mt~\left(\L_\nu^T\right)^{-1} \Big(\nabla_\eta\Phi_\nu\Big)\left(x,\eta_j^\nu\right)
\eeq
 whose Jacobian is non-zero provided that $\Phi$ satisfies the non-degeneracy condition  in (\ref{nondegeneracy}) on the support of $\sigma(x,\L_\nu\eta)$.  

Denote 
$\mathcal{X}=\mathcal{X}(x)\doteq\left(\L_\nu^T\right)^{-1}\Big(\nabla_\eta\Phi_\nu\Big)(x,\eta_j^\nu)$.
By  using (\ref{Omega lj v size}) and (\ref{local diffeo}), we have
\bel{Int Omega lj v size}
\begin{array}{lr}\ds
\int_{\R^\N}\left|\Omega^\nu_{\ell j}(x,y)\right|dx
~\leq~\C_{\sigma~\Phi~N}~2^{-j\left({\N-1\over 2}\right)}2^{\ell\big({m\over 2}\big)}~2^j 2^{j\left({\N-1\over 2}\right)}
\\ \ds~~~~~~~~~~~~~~~~~~~~~~~~~~~~~~~~
\int_{\R^\N} 
\left\{1+ 2^{2j}(\mathcal{X}-y)_1^2+2^j\sum_{i=2}^\N (\mathcal{X}-y)_i^2\right\}^{-N}d\mathcal{X} 
\\\\ \ds~~~~~~~~~~~~~~~~~~~~~~~~~~~~~
~\leq~\C_{\sigma~\Phi~N}~2^{-j\left({\N-1\over 2}\right)}2^{\ell\big({m\over 2}\big)} \iint_{\R\times\R^{\N-1}} \left\{1+ \mathcal{Z}_1^2+     \sum_{i=2}^\N \mathcal{Z}_i^2 \right\}^{-N} 
d\mathcal{Z}_1 \prod_{i=2}^\N d\mathcal{Z}_i 
\\\\ \ds~~~~~~~~~~~~~~~~~~~~~~~~~~~~~~~~~~~~~~~~~~~~~~~~~~~~~~~~~~
\hbox{\small{$ \mathcal{Z}_1=2^{j}(\mathcal{X}-y)_1$,~~ $\mathcal{Z}_i=2^{j/2}(\mathcal{X}-y)_i,~~i=2,\ldots,\N$}}
\\\\ \ds~~~~~~~~~~~~~~~~~~~~~~~~~~~~~
~\leq~\C_{\sigma~\Phi}~2^{-j\left({\N-1\over 2}\right)}2^{\ell\big({m\over 2}\big)}\qquad\hbox{\small{for $N$ sufficiently large.}}
\end{array}
\eeq
Recall {\bf Remark 4.1}. From (\ref{Omega_lj v}) and  (\ref{Int Omega lj v size}), we find
\bel{Int Omega lj size}
\begin{array}{lr}\ds
\int_{\R^\N}\left|\Omega_{\ell j}(x,y)\right|dx~\leq~\sum_{\nu\colon\{\xi^\nu_j\}\subset\mathds{S}^{\N-1}} \int_{\R^\N}\left|\Omega^\nu_{\ell j}(x,y)\right|dx
\\\\ \ds~~~~~~~~~~~~~~~~~~~~~~~~~~~~
~\leq~\C_{\sigma~\Phi}~2^{j\left({\N-1\over 2}\right)} 2^{-\ell m} ~2^{-j\left({\N-1\over 2}\right)} 2^{\ell\big({m\over 2}\big)}
~=~\C_{\sigma~\Phi}~2^{-\ell \big({m\over 2}\big)}.
\end{array}
\eeq
{\bf 2.}~By carrying out the same estimates in (\ref{Theta supp})-(\ref{Int Omega lj size}), we find
\bel{nabla Omega_lj} 
\int_{\R^\N} \left|\nabla_y\Omega_{\ell j}(x,y)\right| dx~\leq~\C_{\sigma~\Phi}~
 2^j2^{-\ell \big({m\over 2}\big)}.
\eeq
This further implies 
\bel{Est2 >} 
\int_{\R^\N} \left|\Omega_{\ell j}(x,y)-\Omega_{\ell j}(x,0)\right| dx~\leq~\C_{\sigma~\Phi}~2^j|y|~2^{-\ell \big({m\over 2}\big)}.
  \eeq

 {\bf 3.}~Recall the region of influence $\Q_{r\ell}$, $0<r<1$ defined in (\ref{R xi natural})-(\ref{Q natural}).

Let $k\in\Z$ and $2^{k-1}< r^{-1}<2^k$. 
  Given $\xi^\nu_j$, $j\ge k$, there exists a $\xi^\mu_k\in \left\{\xi_k^\mu\right\}_\mu\subset\mathds{S}^{\N-1}$ such that $|\xi_j^\nu-\xi_k^\mu|\leq 2^{-k/2}$.

  Suppose $x\in\R^\N\setminus\Q_{r\ell}$, we must have
 \bel{Condition k}
 \left|\left(\nabla_\xi\Phi\right)\left(x,\xi^\mu_k\right)\cdot\xi^\mu_k\right|>\c 2^{-k}\qquad\hbox{or}\qquad 
  \left|\left(\nabla_\xi \Phi\right)\left(x,\xi^\mu_k\right)\right|>\c 2^{-k/2}.
  \eeq 
 Denote $\xi^\mu_k=\L_\nu \eta^\mu_k$. Because $|\eta^\nu_j|=|(\eta^\nu_j)_1|=1$, we find 
 \bel{eta^mu_k}
 1-2^{-k/2}~<~\left|\left(\eta^\mu_k\right)_1\right|~\leq~1,\qquad \left|\left(\eta^\mu_k\right)_i\right|~<~2^{-k/2},\qquad i=2,\ldots,\N.
 \eeq 
 By using the identity in (\ref{xi eta iden}), we have
 \bel{First ineq rewrite}
 \left(\nabla_\xi \Phi\right)\left(x,\xi^\mu_k\right)\cdot\xi^\mu_k~=~\left(\nabla_\eta \Phi_\nu\right) \left(x, \eta^\mu_k\right)\cdot\eta^\mu_k
 \eeq
 where $\Phi_\nu(x,\eta)=\Phi(x,\L_\nu \eta)$ as in (\ref{Phi split}).

The first inequality in (\ref{Condition k}) together with (\ref{eta^mu_k}) imply
 \bel{condition k}
\Big|\left(\partial_{\eta_1}\Phi_\nu\right)\left(x,\eta_k^\mu\right)\Big|>\left[1-{1\over \sqrt{2}}\right]\c 2^{-k}\qquad\hbox{or}\qquad
\sum_{i=2}^\N \Big|\left(\partial_{\eta_i}\Phi_\nu\right)\left(x,\eta_k^\mu\right)\Big|> \c 2^{-k/2}.
\eeq
We have
\bel{nabla Phi diff}
\begin{array}{lr}\ds
\left(\partial_{\eta_1}\Phi_\nu\right)\left(x,\eta_k^\mu\right)-\left(\partial_{\eta_1}\Phi_\nu\right)\left(x,\eta_j^\nu\right)
\\\\ \ds
~=~ \left(\nabla_\eta\partial_{\eta_1}\Phi_\nu\right)\left(x,\eta^\nu_j\right)\cdot\left( \eta_k^\mu- \eta_j^\nu\right)+\O(1)  \left|\eta_k^\mu- \eta_j^\nu\right|^2
 \\\\ \ds
 ~=~ \left(\partial_{\eta_1}\nabla_\eta\Phi_\nu\right)\left(x, \eta^\nu_j\right)\cdot\left( \eta_k^\mu- \eta_j^\nu\right)+\O(1)  \left|\eta_k^\mu- \eta_j^\nu\right|^2. 
\end{array}
\eeq
Note that $\nabla_\eta\Phi_\nu$ is homogeneous of degree $0$ in $\eta$. 
By using the identity of $\partial_{\eta_1}$ in (\ref{radial derivative}), we find
\bel{condition1}
\begin{array}{lr}\ds
\left|\left(\partial_{\eta_1}\Phi_\nu\right)\left(x,\eta_k^\mu\right)-\left(\partial_{\eta_1}\Phi_\nu\right)\left(x,\eta_j^\nu\right)\right|~\leq~\C ~\left(2^{-k/2} \right)^2
~=~\C ~2^{-k}.
\end{array}
\eeq
On the other hand, the mean value theorem implies
\bel{condition2}
\left| \left(\nabla_\eta\Phi_\nu\right)\left(x,\eta_k^\mu\right)-\left(\nabla_\eta\Phi_\nu\right)\left(x,\eta_j^\nu\right)\right|~\leq~\C~ 2^{-k/2}.
\eeq
Suppose $y\in\B_r, 0<r<1$. 
From (\ref{condition k}), (\ref{condition1}) and (\ref{condition2}), by using the triangle inequality and taking into account for $\left|\L_\nu^T y\right|\leq 2^{-k}$, we obtain
\bel{condition jk}
2^{2j}\left[\left(\nabla_\eta\Phi_\nu\right)\left(x,\eta_j^\nu\right)-\L_\nu^T y\right]_1^2+4\pi^2 2^{j}\sum_{i=2}^\N \left[\left(\nabla_\eta\Phi_\nu\right)\left(x,\eta_j^\nu\right)-\L_\nu^T y\right]_i^2~\ge~2^{j-k}
\eeq
provided that $\c>0$ inside (\ref{condition k}) is large enough.

Furthermore, the second inequality in (\ref{Condition k}) and (\ref{xi eta iden}) imply
$\left| \left(\nabla_\eta\Phi_\nu\right) \left(x,\eta_k^\mu\right)\right|>\c2^{-k/2}$.
Together with (\ref{condition2}), by using the triangle inequality and taking into account for $\left|\L_\nu^T y\right|\leq 2^{-k}$, we find (\ref{condition jk}) again if $\c>0$ is sufficiently large.

Consider $j\ge k$ and $x\in\R^\N\setminus\Q_{r\ell}$. 
From (\ref{Omega lj v size}) and (\ref{condition jk}), we have 
\bel{Omega lj v size j-k}
\begin{array}{lr}\ds
\left|\Omega^\nu_{\ell j}(x,y)\right|~\leq~\C_{N}~2^{-j\left({\N-1\over 2}\right)}2^{\ell\big({m\over 2}\big)}~2^j 2^{j\left({\N-1\over 2}\right)} 
\\ \ds~~~~~~~~~~~~~~~~~~~~
\left\{1+4\pi^2 2^{2j}\left[\left(\nabla_\eta\Phi_\nu\right)\left(x,\eta_j^\nu\right)-\L_\nu^T y\right]_1^2+4\pi^2 2^{j}\sum_{i=2}^\N \left[\left(\nabla_\eta\Phi_\nu\right)\left(x,\eta_j^\nu\right)-\L_\nu^T y\right]_i^2 \right\}^{-N}
\\\\ \ds~~~~~~~~~~~~~~~~~
~\leq~\C_{N}~2^{-j\left({\N-1\over 2}\right)}2^{\ell\big({m\over 2}\big)}~2^j 2^{j\left({\N-1\over 2}\right)} ~2^{-j+k}
\\ \ds~~~~~~~~~~~~~~~~~~~~
\left\{1+4\pi^2 2^{2j}\left[\left(\nabla_\eta\Phi_\nu\right)\left(x,\eta_j^\nu\right)-\L_\nu^T y\right]_1^2+4\pi^2 2^{j}\sum_{i=2}^\N \left[\left(\nabla_\eta\Phi_\nu\right)\left(x,\L_\nu\eta_j^\nu\right)-\L_\nu^T y\right]_i^2 \right\}^{1-N}.
\end{array}
\eeq
Now,  repeat the estimates from  (\ref{Theta supp}) to (\ref{Int Omega lj size})
expect that (\ref{Omega lj v size}) is replaced by (\ref{Omega lj v size j-k}). Finally, we arrive at
\bel{Result Q l}
\int_{\R^\N\setminus \Q_{r\ell}} \left|\Omega_{\ell j}(x,y) \right|dx~\leq~\C~
{2^{-j}\over r} ~ 2^{-\ell \big({m\over 2}\big)}.
\eeq

{\bf Case Two}~~Consider $j/2\leq\ell<j$. Write $\xi=(\tau,\lambda)\in\R^n\times\R^m$. Without lose of the generality,  take $\tau=(\xi_1,\ldots,\xi_n)$ and $\lambda=(\xi_{n+1},\ldots,\xi_{n+m})$. $\left\{ \tau^\nu_j\right\}_\nu$ is a collection of points  equally distributed on  $\mathds{S}^{n-1}$ with a grid length between $2^{-j/2-1}$ and $2^{-j/2}$. Denote $\xi^\nu_j=(\tau^\nu_j,0)\in\R^n\times\R^m$.

Consider
$\xi=\L_\nu \eta$ where $\L_\nu$ is an $\N\times \N$-orthogonal matrix with $\det\L_\nu=1$. 
Furthermore, we have $\eta^1=(\eta_1,0)\in\R\times\R^{\N-1}$ in the same direction of $\xi^\nu_j$. Write $\eta=(\zeta,\gamma)\in\R^n\times\R^m$ for $\zeta=(\eta_1,\ldots,\eta_n)$ and $\gamma=(\eta_{n+1},\ldots,\eta_{n+m})$.
We further require
\bel{L_nu block}
\L_\nu~=~\left[\begin{array}{lr}\ds L_\nu 
\\ \ds
~~~~~ Id\end{array}\right],\qquad \tau~=~L_\nu \zeta,\qquad \lambda~=~\gamma
\eeq
where $L_\nu$ is an $n\times n$-orthogonal matrix with $\det L_\nu=1$.

Recall (\ref{Omega_j sum+}). We write
\bel{Omega_lj v+}
\begin{array}{cc}\ds
\Omega_{\ell j}(x,y)~=~\sum_{\nu\colon\{ \tau^\nu_j\}_\nu\subset\mathds{S}^{n-1}} {^+}\Omega^\nu_{\ell j}(x,y),
\\\\ \ds
{^+}\Omega^\nu_{\ell j}(x,y)~=~\int_{\R^\N} e^{2\pi\i \big[\Phi(x,\L_\nu\eta)-y\cdot\L_\nu\eta \big]}\delta_\ell(\L_\nu\eta)\phi_j(\L_\nu\eta){^+}\varphi^\nu_j(\L_\nu\eta)\sigma(x,\L_\nu\eta)d\eta
\end{array}
\eeq
and
\bel{Omega_j v+}
\begin{array}{cc}\ds
{^+}\Omega_j(x,y)~=~\sum_{\nu\colon\{ \tau^\nu_j\}_\nu\subset\mathds{S}^{n-1}} {^+}\Omega^\nu_j(x,y),
\\\\ \ds
{^+}\Omega^\nu_j(x,y)~=~\int_{\R^\N} e^{2\pi\i \big[\Phi(x,\L_\nu\eta)-y\cdot\L_\nu\eta \big]}\mathfrak{S}^+_j(\L_\nu\eta)\phi_j(\L_\nu\eta){^+}\varphi^\nu_j(\L_\nu\eta)\sigma(x,\L_\nu\eta)d\eta
\end{array}
\eeq
where $\varphi^\nu_j$ and $\mathfrak{S}^+_j$ are defined in (\ref{phi^v_j +}) and (\ref{delta_j+-}). 

Denote $\eta^\nu_j=(1,0)\in\R\times\R^{\N-1}$ so that $\xi^\nu_j=\L_\nu \eta^\nu_j$. 
Similar to (\ref{Theta}), we write
\bel{Theta+lj}
\begin{array}{cc}\ds
{^+}\Omega^\nu_{\ell j}(x,y)~=~\int_{\R^\N} e^{2\pi\i\Big[\left(\nabla_\eta\Phi_\nu\right)\big(x,\eta_j^\nu\big)-\L_\nu^T y\Big]\cdot\eta}{^+}\Theta^\nu_{\ell j}(x,\eta)d\eta,
\\\\ \ds
{^+}\Theta^\nu_{\ell j}(x,\eta)~=~e^{2\pi\i\Psi_\nu(x,\eta)}\delta_\ell(\L_\nu\eta)\phi_j(\L_\nu\eta){^+}\varphi^\nu_j(\L_\nu\eta)\sigma(x,\L_\nu\eta)
\end{array}
\eeq
and
\bel{Theta+j}
\begin{array}{cc}\ds
{^+}\Omega^\nu_j(x,y)~=~\int_{\R^\N} e^{2\pi\i\Big[\left(\nabla_\eta\Phi_\nu\right)\big(x,\eta_j^\nu\big)-\L_\nu^T y\Big]\cdot\eta}{^+}\Theta^\nu_j(x,\eta)d\eta,
\\\\ \ds
{^+}\Theta^\nu_j(x,\eta)~=~e^{2\pi\i\Psi_\nu(x,\eta)}\mathfrak{S}^+_j(\L_\nu\eta)\phi_j(\L_\nu\eta){^+}\varphi^\nu_j(\L_\nu\eta)\sigma(x,\L_\nu\eta)
\end{array}
\eeq
where $\Psi_\nu$ is defined in (\ref{Phi split}).

Define 
\bel{D+lj}
{^+}\D_{\ell j}~=~I+2^{2j} \partial_{\eta_1}^2+2^j  \sum_{i=2}^n \partial^2_{\eta_i}+2^{2(j-\ell)}\sum_{k=n+1}^{n+m} \partial^2_{\eta_k},\qquad j>0,~~j/2\leq\ell<j
\eeq
and
\bel{D+j}
{^+}\D_j~=~I+2^{2j} \partial_{\eta_1}^2+2^j  \sum_{i=2}^n \partial^2_{\eta_i}+\sum_{k=n+1}^{n+m} \partial^2_{\eta_k},\qquad j>0.
\eeq
\begin{lemma} Suppose $\sigma\in\S^{-{\N-1\over 2}}$.
Given $j>0$ and $j/2\leq\ell<j$, we have
\bel{Theta Diff Ineq+lj}
\left|{^+}\D_{\ell j}^N {^+}\Theta^\nu_{\ell j}(x,\eta)\right|~\leq~\C_N~2^{-j\left({\N-1\over 2}\right)}2^{\ell\big[{m\over 2}-\ve\big]},\qquad N\ge0
\eeq
and
\bel{Theta Diff Ineq+j}
\left|{^+}\D_j^N {^+}\Theta^\nu_j(x,\eta)\right|~\leq~\C_N~2^{-j\left[{n-1\over 2}+\ve\right]},\qquad N\ge0
\eeq
for some $\ve=\ve(\rho_1,\rho_2)>0$.
\end{lemma}

{\bf Proof}~~First, ${^+}\varphi^\nu_j$ defined in (\ref{phi^v_j +}) has a support  inside $\Gamma^\nu_j$. The identity of $\partial_{\eta_1}$ in (\ref{radial derivative}) remains valid.  Consequently,  ${^+}\varphi^\nu_j$
also satisfies the estimates in (\ref{partial 0})-(\ref{nu Diff est 1}). We have
\bel{var^v_j+ Diff Ineq}
\left|\partial^\alpha_\eta\left[{^+}\varphi^\nu_j(\L_\nu \eta)\right]\right|~\leq~\C_\alpha~2^{(j/2)\big[|\alpha|-\alpha_1\big]} |\eta|^{-|\alpha|}
\eeq
for every multi-index $\alpha$.

Recall $\delta_\ell$ satisfying the differential inequality in (\ref{delta_l Diff Ineq}). Moreover,
$\supp\delta_\ell\subset\Lambda_\ell$ which is defined in (\ref{Cone}).  Write $\xi=(\tau,\lambda)\in\R^n\times\R^m$.    
We find $|\tau|\sim2^j$,  $|\lambda|\sim2^{j-\ell}$ whenever $\xi\in \Lambda_\ell\cap\left\{2^{j-1}\leq|\xi|\leq2^{j+1}\right\}$. 

On the other hand, $\mathfrak{S}_j^+$ satisfies the first differential inequality in (\ref{delta_j Diff Ineq}). Furthermore, $\supp\mathfrak{S}_j^+$ is contained inside the region where $|\lambda|\leq 2^{-j+1}|\tau|$.

Suppose  that $\sigma\in\S^{-{\N-1\over2}}$ satisfies the differential inequality in (\ref{Class})-(\ref{rho_1, rho_2}). We find
 $\left({1\over1+|\tau|}\right)^{\rho_1}\left({1\over 1+|\lambda|}\right)^{\rho_2}=\left({1\over1+|\tau|}\right)^{{n-1\over 2}+\ve}\left({1\over 1+|\lambda|}\right)^{{m\over 2}-\ve}$  for $0<\ve<{m\over 2}$ where $\rho_1= {n-1\over 2}+\ve$, $\rho_2>{m-1\over 2}$ and  $\rho_1+\rho_2={\N-1\over 2}$.
All together, we have
\bel{Ineq+lj}
\begin{array}{lr}\ds
\left| \partial^\alpha_\tau\partial^\beta_\lambda \left[\delta_\ell(\xi) \phi_j(\xi) \sigma(x,\xi)\right]\right|
~\leq~\C_{\alpha~\beta} \left({1\over1+|\tau|}\right)^{{n-1\over 2}+\ve}\left({1\over 1+|\lambda|}\right)^{{m\over 2}-\ve} |\tau|^{-|\alpha|}|\lambda|^{-|\beta|}
\\\\ \ds~~~~~~~~~~~~~~~~~~~~~~~~~~~~~~~~~~~~~~~~~~
~\leq~\C_{\alpha~\beta}~2^{-j\big({\N-1\over 2}\big)}2^{\ell\big[{m\over 2}-\ve\big]} 2^{-j |\alpha|} 2^{-(j-\ell)|\beta|},
\\\\ \ds
\left| \partial^\alpha_\tau\partial^\beta_\lambda \left[\mathfrak{S}_j^+(\xi) \phi_j(\xi) \sigma(x,\xi)\right]\right|
~\leq~\C_{\alpha~\beta} \left({1\over1+|\tau|}\right)^{{n-1\over 2}+\ve}\left({1\over 1+|\lambda|}\right)^{{m\over 2}-\ve} |\tau|^{-|\alpha|}
\\\\ \ds~~~~~~~~~~~~~~~~~~~~~~~~~~~~~~~~~~~~~~~~~~~~
~\leq~\C_{\alpha~\beta}~2^{-j\big[{n-1\over 2}+\ve\big]}2^{-j |\alpha|} 
\end{array}
\eeq
for every multi-indices $\alpha,\beta$.

Let $\L_\nu$ defined in (\ref{L_nu block}). Write $\eta=(\zeta,\gamma)$. We have $\tau= L_\nu \zeta$ and $\lambda=\gamma$. 

From (\ref{var^v_j+ Diff Ineq}) and (\ref{Ineq+lj}), a direct computation shows
\bel{INEQ+lj}
\begin{array}{lr}
\left|\p_\zeta^\alpha\p_\gamma^\beta \left[\delta_\ell (\L_\nu \eta)\phi_j(\L_\nu\eta){^+}\varphi_j^\nu \left(\L_\nu\eta\right)\sigma(x,\L_\nu \eta)\right] \right|
\\\\ \ds
~\leq~\C_{\alpha~\beta}~2^{-j\left({\N-1\over 2}\right)}2^{\ell\big[{m\over 2}-\ve\big]} 2^{-j \alpha_1} 2^{-(j/2)\big[|\alpha|-\alpha_1\big]}2^{-(j-\ell)|\beta|},
\\\\\\ \ds
\left|\p_\zeta^\alpha\p_\gamma^\beta \left[\mathfrak{S}_j^+(\L_\nu\eta) (\L_\nu \eta)\phi_j(\L_\nu\eta){^+}\varphi_j^\nu \left(\L_\nu\eta\right)\sigma(x,\L_\nu \eta)\right] \right|
\\\\ \ds
~\leq~\C_{\alpha~\beta}~2^{-j\big[{n-1\over 2}+\ve\big]} 2^{-j \alpha_1} 2^{-(j/2)\big[|\alpha|-\alpha_1\big]}
\end{array}
\eeq
for every multi-indices $\alpha,\beta$.

Recall $\Psi_\nu$ satisfying the differential inequality in (\ref{d Est Psi}). Together with the estimates in (\ref{INEQ+lj}), we obtain (\ref{Theta Diff Ineq+lj})-(\ref{Theta Diff Ineq+j}) as desired. 
\endproof

Let ${^+}\Omega^\nu_{\ell j}(x,y)$ and ${^+}\Theta^\nu_{\ell j}(x,y)$ for $j/2\leq\ell<j$ defined in (\ref{Theta+lj}). 
We find 
\bel{supp Theta+}
\left|\supp {^+}\Theta^\nu_{\ell j}(x,\eta)\right|~\leq~\C~2^j 2^{j\left({n-1\over 2}\right)} 2^{(j-\ell)m}.
\eeq
By applying (\ref{Theta Diff Ineq+lj}), an $N$-fold integration by parts associated to ${^+}\D_{\ell j}$ defined in (\ref{D+lj}) shows 
\bel{Omega+lj v size}
\begin{array}{lr}\ds
\left|{^+}\Omega^\nu_{\ell j}(x,y)\right|~\leq~
\left\{1+4\pi^2 2^{2j}\left[\left(\nabla_\eta\Phi_\nu\right)\left(x,\eta_j^\nu\right)-\L_\nu^T y\right]_1^2
+
4\pi^2 2^j\sum_{i=2}^n \left[\left(\nabla_\eta\Phi_\nu\right)\left(x,\eta_j^\nu\right)-\L_\nu^T y\right]_i^2\right.
\\ \ds~~~~~~~~~~~~~~~~~~~~~~~~~
\left.~+4\pi^2 2^{2(j-\ell)} \sum_{k=n+1}^{n+m} \left[\left(\nabla_\eta\Phi_\nu\right)\left(x,\eta_j^\nu\right)-\L_\nu^T y\right]_k^2\right\}^{-N}\int_{\R^\N} \left| {^+}\D_{\ell j}^N {^+}\Theta^\nu_{\ell j}(x,\eta)\right|d\eta
\\\\ \ds~~~~~~~~~~~~~~~~~~
~\leq~\C_{N}~2^{-j\left({\N-1\over 2}\right)}2^{\ell\big[{m\over 2}-\ve\big]}~2^j 2^{j\left({n-1\over 2}\right)} 2^{(j-\ell)m}
\\ \ds~~~~~~~~~~~~~~~~~~~~~~~~~~
\left\{1+4\pi^2 2^{2j}\left[\left(\nabla_\eta\Phi_\nu\right)\left(x,\eta_j^\nu\right)-\L_\nu^T y\right]_1^2
+
4\pi^2 2^j\sum_{i=2}^n \left[\left(\nabla_\eta\Phi_\nu\right)\left(x,\eta_j^\nu\right)-\L_\nu^T y\right]_i^2\right.
\\ \ds~~~~~~~~~~~~~~~~~~~~~~~~~
\left.~+4\pi^2 2^{2(j-\ell)} \sum_{k=n+1}^{n+m} \left[\left(\nabla_\eta\Phi_\nu\right)\left(x,\eta_j^\nu\right)-\L_\nu^T y\right]_k^2\right\}^{-N}
.
\end{array}
\eeq
Recall $\mathcal{X}_\Phi\colon x\mt\left(\L_\nu^T\right)^{-1}\left(\nabla_\eta\Phi_\nu\right)\left(x,\eta_j^\nu\right)$  from (\ref{local diffeo}) which is a local diffeomorphism 
whose Jacobian is non-zero because $\Phi$ satisfies (\ref{nondegeneracy}) on the support of $\sigma(x,\L_\nu \eta)$. 

Denote
$\mathcal{X}=\mathcal{X}(x)=\left(\L_\nu^T\right)^{-1}\left(\nabla_\eta\Phi_\nu\right)\left(x,\eta_j^\nu\right)$.
By  using (\ref{Omega+lj v size}), we have
\bel{Int Omega+lj v size}
\begin{array}{lr}\ds
\int_{\R^\N}\left|{^+}\Omega^\nu_{\ell j}(x,y)\right|dx
~\leq~
\C_N~2^{-j\left({\N-1\over 2}\right)}2^{\ell\big[{m\over 2}-\ve\big]}~2^j 2^{j\left({n-1\over 2}\right)} 2^{(j-\ell)m}
\\ \ds~~~~~~~~~~~~~~~~~~~~~~~~~~~~~~
\iint_{\R^n\times\R^m} 
\left\{1+ 2^{2j}(\mathcal{X}-y)_1^2+2^j\sum_{i=2}^n (\mathcal{X}-y)_i^2+2^{2(j-\ell)}\sum_{k=n+1}^{n+m} (\mathcal{X}-y)_i^2\right\}^{-N}
d\mathcal{X}
\\\\ \ds
~\leq~\C_{N}~2^{-j\left({\N-1\over 2}\right)}2^{\ell\big[{m\over 2}-\ve\big]} \iint_{\R^n\times\R^m} \left\{1+ \mathcal{Z}_1^2+     \sum_{i=2}^n \mathcal{Z}_i^2 + \sum_{k=n+1}^{n+m} \mathcal{Z}_i^2\right\}^{-N} 
d\mathcal{Z}_1 \cdots d\mathcal{Z}_{n+m} 
\\\\ \ds~~~~~
\hbox{\small{$ \mathcal{Z}_1=2^{j}(\mathcal{X}-y)_1$,~~~~ $\mathcal{Z}_i=2^{j/2}(\mathcal{X}-y)_i,~i=2,\ldots,n$,~~~~$\mathcal{Z}_k=2^{j-\ell}(\mathcal{X}-y)_k,~k=n+1,\ldots,n+m$}}
\\\\ \ds
~\leq~\C~2^{-j\left({\N-1\over 2}\right)}2^{\ell\big[{m\over 2}-\ve\big]}\qquad\hbox{\small{for $N$ sufficiently large.}}
\end{array}
\eeq
Recall {\bf Remark 4.3}. From (\ref{Omega_lj v+}) and (\ref{Int Omega+lj v size}), we find
\bel{Int Omega+lj size}
\begin{array}{lr}\ds
\int_{\R^\N}\left|\Omega_{\ell j}(x,y)\right|dx~\leq~\sum_{\nu\colon\{\tau^\nu_j\}\subset\mathds{S}^{n-1}} \int_{\R^\N}\left|{^+}\Omega^\nu_{\ell j}(x,y)\right|dx
\\\\ \ds~~~~~~~~~~~~~~~~~~~~~~~~~~~~
~\leq~\C~2^{j\left({n-1\over 2}\right)}  ~2^{-j\left({\N-1\over 2}\right)} 2^{\ell\big[{m\over 2}-\ve\big]}
~=~\C~2^{-(j-\ell) \big[{m\over 2}-\ve\big]}~2^{-j\ve}.
\end{array}
\eeq

Let ${^+}\Omega^\nu_j(x,y)$ and ${^+}\Theta^\nu_j(x,y)$ defined in (\ref{Theta+j}). 
We find 
\bel{supp Theta+j}
\left|\supp {^+}\Theta^\nu_j(x,\eta)\right|~\leq~\C~2^j 2^{j\left({n-1\over 2}\right)}.
\eeq
By applying (\ref{Theta Diff Ineq+j}), an $N$-fold integration by parts associated to ${^+}\D_j$ defined in (\ref{D+j}) shows 
\bel{Omega+j v size}
\begin{array}{lr}\ds
\left|{^+}\Omega^\nu_j(x,y)\right|~\leq~
\left\{1+4\pi^2 2^{2j}\left[\left(\nabla_\eta\Phi_\nu\right)\left(x,\eta_j^\nu\right)-\L_\nu^T y\right]_1^2
+
4\pi^2 2^j\sum_{i=2}^n \left[\left(\nabla_\eta\Phi_\nu\right)\left(x,\eta_j^\nu\right)-\L_\nu^T y\right]_i^2\right.
\\ \ds~~~~~~~~~~~~~~~~~~~~~~~~~~~
\left.+4\pi^2 \sum_{k=n+1}^{n+m} \left[\left(\nabla_\eta\Phi_\nu\right)\left(x,\eta_j^\nu\right)-\L_\nu^T y\right]_k^2\right\}^{-N}\int_{\R^\N} \left| {^+}\D_j^N {^+}\Theta^\nu_j(x,\eta)\right|d\eta
\\\\ \ds~~~~~~~~~~~~~~~~~~
~\leq~\C_{N}~2^{-j\left[{n-1\over 2}+\ve\right]}~2^j 2^{j\left({n-1\over 2}\right)} 
\\ \ds~~~~~~~~~~~~~~~~~~~~~~~~~~
\left\{1+4\pi^2 2^{2j}\left[\left(\nabla_\eta\Phi_\nu\right)\left(x,\eta_j^\nu\right)-\L_\nu^T y\right]_1^2
+
4\pi^2 2^j\sum_{i=2}^n \left[\left(\nabla_\eta\Phi_\nu\right)\left(x,\eta_j^\nu\right)-\L_\nu^T y\right]_i^2\right.
\\ \ds~~~~~~~~~~~~~~~~~~~~~~~~~~
\left.+4\pi^2 \sum_{k=n+1}^{n+m} \left[\left(\nabla_\eta\Phi_\nu\right)\left(x,\eta_j^\nu\right)-\L_\nu^T y\right]_k^2\right\}^{-N}.
\end{array}
\eeq
Denote
$\mathcal{X}=\mathcal{X}(x)=\left(\L_\nu^T\right)^{-1}\left(\nabla_\eta\Phi_\nu\right)\left(x,\eta_j^\nu\right)$.
By  using (\ref{Omega+j v size}), we have
\bel{Int Omega+j v size}
\begin{array}{lr}\ds
\int_{\R^\N}\left|{^+}\Omega^\nu_j(x,y)\right|dx
~\leq~
\C_{N}~2^{-j\left[{n-1\over 2}+\ve\right]}~2^j 2^{j\left({n-1\over 2}\right)} 
\\ \ds~~~~~~~~~~~~~~~~~~~~~~~~~~~~~~~~~~~~
\iint_{\R^n\times\R^m} 
\left\{1+ 2^{2j}(\mathcal{X}-y)_1^2+2^j\sum_{i=2}^n (\mathcal{X}-y)_i^2+\sum_{k=n+1}^{n+m} (\mathcal{X}-y)_k^2\right\}^{-N}
d\mathcal{X}
\\\\ \ds
~\leq~\C_{N}~2^{-j\left[{n-1\over 2}+\ve\right]}\iint_{\R^n\times\R^m} \left\{1+ \mathcal{Z}_1^2+     \sum_{i=2}^n \mathcal{Z}_i^2 +\sum_{k=n+1}^{n+m} \mathcal{Z}_k^2\right\}^{-N} 
d\mathcal{Z}_1 \cdots d\mathcal{Z}_{n+m}
\\\\ \ds~~~~~~~
\hbox{\small{$ \mathcal{Z}_1=2^{j}(\mathcal{X}-y)_1$,~~~~ $\mathcal{Z}_i=2^{j/2}(\mathcal{X}-y)_i,~i=2,\ldots,n$,~~~~ $\mathcal{Z}_k=(\mathcal{X}-y)_k,~k=n+1,\ldots,n+m$}}
\\\\ \ds
~\leq~\C~2^{-j\left[{n-1\over 2}+\ve\right]}\qquad\hbox{\small{for $N$ sufficiently large.}}
\end{array}
\eeq
Recall {\bf Remark 4.3}. From (\ref{Omega_j v+}) and (\ref{Int Omega+j v size}), we find
\bel{Int Omega+j size}
\begin{array}{lr}\ds
\int_{\R^\N}\left|{^+}\Omega_j(x,y)\right|dx~\leq~\sum_{\nu\colon\{\tau^\nu_j\}\subset\mathds{S}^{n-1}} \int_{\R^\N}\left|{^+}\Omega^\nu_j(x,y)\right|dx
\\\\ \ds~~~~~~~~~~~~~~~~~~~~~~~~~~~~~
~\leq~\C~2^{j\left({n-1\over 2}\right)}  ~2^{-j\left[{n-1\over 2}+\ve\right]} 
~=~\C~2^{-j\ve}.
\end{array}
\eeq
Lastly, recall ${^\sharp}\Omega_j$ defined in (\ref{F sharp}). By using (\ref{Int Omega+lj size}) and (\ref{Int Omega+j size}), we have
\bel{sharp Omega_j EST} 
\begin{array}{lr}\ds
\int_{\R^\N}\left|{^\sharp}\Omega_j(x,y)\right|dx~\leq~\sum_{j/2\leq\ell<j}\int_{\R^\N} \left|\Omega_{\ell j}(x,y)\right|dx+\int_{\R^\N}\left|{^+}\Omega_j(x,y)\right|dx
\\\\ \ds~~~~~~~~~~~~~~~~~~~~~~~~~~~~~
~\leq~\C \sum_{j/2\leq\ell<j} ~2^{-(j-\ell) \big[{m\over 2}-\ve\big]}~2^{-j\ve}+ 2^{-j \ve}
\\\\ \ds~~~~~~~~~~~~~~~~~~~~~~~~~~~~~
~\leq~\C~2^{-j\ve/2},\qquad \hbox{\small{$0<\ve<{m\over 2}$}}.
\end{array}
\eeq

\section{$\H^1\mt\L^1$-boundedness of the adjoint operator}
\setcounter{equation}{0}
Let $\F$ defined in (\ref{Ff})-(\ref{nondegeneracy}). We find
\bel{F*f}
\F^* f(x)~=~\int_{\R^\N} f(y) \left\{ \int_{\R^\N} e^{2\pi\i\big[x\cdot\xi-\Phi(y,\xi)\big]} \bar{\sigma}(y,\xi) d\xi\right\} dy.
\eeq
Observe that $\F^*$ does not have a compact support in $x$. 

Let $a$ be an $\H^1$-atom. Without lose of the generality, assume $\supp a\subset\B_r$. 
We claim
\bel{H^1 to L^1 EST*}
\int_{\R^\N} \left| \F^* f(x)\right|dx~\leq~\C_{\rho_1~\rho_2}~\left\| f\right\|_{\H^1(\R^\N)},
\qquad\sigma\in\S^{-{\N-1\over2}}.
\eeq
Let $\F_o$ and ${^\natural}\F_\ell, {^\sharp}\F_j, {^\flat}\F_j$  defined in (\ref{F_o}) and (\ref{F natural})-(\ref{F flat}) respectively.  Recall {\bf Remark 2.2}. We have
\bel{F* SUM}
\begin{array}{cc}\ds
\F^* f(x)~=~\F^*_o f(x)+{^\natural}\F^* f(x)+{^\sharp}\F^* f(x)+{^\flat}\F^* f(x).
\end{array}
\eeq
In particular, 
\bel{F*_o}
\begin{array}{cc}\ds
\F^*_o f(x)~=~\int_{\R^\N} f(y)\Omega^*_o(x,y)dy,
\\\\ \ds
\Omega^*_o(x,y)~=~
\int_{\R^\N} e^{2\pi\i \big[x\cdot\xi-\Phi(y,\xi) \big]}\bar{\phi}(\xi)\bar{\sigma}(y,\xi)d\xi
\end{array}
\eeq
and
\bel{F*_l natural}
\begin{array}{cc}\ds
{^\natural}\F^* f(x)~=~\sum_{j>0}\sum_{\ell\colon-j/2<\ell<j/2} \F^*_{\ell j} f(x),
\\\\ \ds
\F^*_{\ell j} f(x)~=~\int_{\R^\N} f(y)\Omega^*_{\ell j}(x,y)dy,
\\\\ \ds
\Omega^*_{\ell j}(x,y)~=~
\int_{\R^\N} e^{2\pi\i \big[x\cdot\xi-\Phi(y,\xi) \big]}\bar{\delta}_\ell(\xi)\bar{\phi}_j(\xi)\bar{\sigma}(y,\xi)d\xi
\end{array}
\eeq
where $\delta_\ell$, $\phi$ and $\phi_j$ are defined in (\ref{delta_t}) and (\ref{phi_j}).

Let $\ds\U=\left\{\nabla_\xi\Phi(y,\xi)\colon y\in\supp\sigma\right\}$ where $\nabla_\xi\Phi$ is homogeneous of degree $0$ in $\xi$. 
Define
\bel{compact set}
\mathfrak{Q}~=~\left\{x\in\R^\N\colon \dist\left\{ x,\mathcal{U}\right\}\leq1\right\}.
\eeq
Suppose $x\in\R^\N\setminus\mathfrak{Q}$. We must have
$\left|x-\nabla_\xi\Phi(y,\xi)\right|\sim|x|$ as $|x|\mt\infty$ whenever $y\in\supp\bar{\sigma}$.

Let $\Omega^*_o$ be defined in (\ref{F*_o}).  An $N$-fold integration by parts shows
\bel{Kernel int by parts Omega_o}
\begin{array}{lr}\ds
\Omega^*_o(x,y)~=~
(2\pi\i)^{-N} \int_{\R^\N} e^{2\pi\i\big[x\cdot\xi-\Phi(y,\xi)\big]} \nabla_\xi^N\Bigg\{ \left[ x-\nabla_\xi\Phi(y,\xi)\right]^{-N} \bar{\phi}(\xi)\bar{\sigma}(y,\xi) \Bigg\}d\xi.
\end{array}
\eeq
Note that  $\supp\phi$ is contained inside the ball: $|\xi|\leq2$. We find $|\Omega^*_o(x,y)|$ uniformly bounded.
Moreover, the integrand inside (\ref{Kernel int by parts Omega_o}) has norm bounded by $\C_N\left({1\over 1+|x|}\right)^N$ whenever $x\in\R^\N\setminus \mathfrak{Q}$.
Therefore, we have
\bel{F*_o L^1}
\left\| \F^*_o a\right\|_{\L^1(\mathfrak{Q})}+\left\| \F^*_o a\right\|_{\L^1(\R^\N\setminus\mathfrak{Q})}~\leq~\C~\left\| a\right\|_{\L^1(\R^\N)}.
\eeq
Let $\Omega^*_{\ell j}$ be defined in (\ref{F*_l natural}) for $j>0, -j/2\ell<j/2$.  An $N$-fold integration by parts shows
\bel{Kernel int by parts Omega_l}
\begin{array}{lr}\ds
\Omega^*_{\ell j}(x,y)~=~
(2\pi\i)^{-N} \int_{\R^\N} e^{2\pi\i\big[x\cdot\xi-\Phi(y,\xi)\big]} \nabla_\xi^N\Bigg\{ \left[ x-\nabla_\xi\Phi(y,\xi)\right]^{-N} \bar{\delta}_\ell(\xi)\bar{\phi}_j(\xi)\bar{\sigma}(y,\xi) \Bigg\}d\xi.
\end{array}
\eeq
Because $\nabla_\xi\Phi(x,\xi)$ is homogeneous of degree $0$ in $\xi$, we have $\left|\nabla_\xi^N \left[ \nabla_\xi\Phi(x,\xi)\right]\right|\leq\C_N |\xi|^{-N}, N\ge0$ as discussed earlier in section 5. Moreover, $\bar{\delta}_\ell$ satisfies the differential inequality in (\ref{delta_l Diff Ineq}).  $\sigma\in\S^{-{\N-1\over2}}$ satisfies the differential inequality in (\ref{Class})-(\ref{rho_1, rho_2}). For $\xi\in\supp\bar{\delta}_\ell\cap\supp\bar{\phi}_j$, we find $|\tau|\sim2^j$,  $|\lambda|\sim2^{j-\ell}$ for $\ell\ge0$ and $|\tau|\sim2^{j+\ell}$,  $|\lambda|\sim2^{j}$ for $\ell\leq0$.

From direct computation, we have
\bel{Omega*_lj norm>}
\begin{array}{lr}\ds
\left| \Omega^*_{\ell j}(x,y)\right|~\leq~\C_N~\left({1\over 1+|x|}\right)^N~2^{(j-\ell)N} ~2^{jn}2^{(j-\ell)m}
\\\\ \ds~~~~~~~~~~~~~~~~~
~\leq~\C_N~\left({1\over 1+|x|}\right)^N~2^{jN/2} ~2^{j\N}2^{-\ell m}\qquad\hbox{\small{ if ~$0\leq\ell<j/2$}}
\end{array}
\eeq
and
\bel{Omega*_lj norm<}
\begin{array}{lr}\ds
\left| \Omega^*_{\ell j}(x,y)\right|~\leq~\C_N~\left({1\over 1+|x|}\right)^N~2^{(j+\ell)N} ~2^{(j+\ell)n}2^{jm}
\\\\ \ds~~~~~~~~~~~~~~~~~
~\leq~\C_N~\left({1\over 1+|x|}\right)^N~2^{jN/2} ~2^{j\N}2^{-\ell n}\qquad\hbox{\small{ if ~$-j/2<\ell\leq0$}}
\end{array}
\eeq
whenever $x\in\R^\N\setminus\mathfrak{Q}$.

By taking $N=2\N+1$ in (\ref{Omega*_lj norm>})-(\ref{Omega*_lj norm<}), we find 
\bel{Omega*_lj norm}
\left| \Omega^*_{\ell j}(x,y)\right|~\leq~\C~\left({1\over 1+|x|}\right)^{2\N+1} ~2^{-j/2} \left\{\begin{array}{lr}\ds 2^{-m\ell},\qquad 0\leq\ell<j/2,
\\ \ds
2^{n\ell},\qquad -j/2<\ell\leq0,
\end{array}\right.
\qquad \hbox{\small{$x\in\R^\N\setminus\mathfrak{Q}$}}.
\eeq
Recall ${^\natural}\F$ defined in (\ref{F*_l natural}). Clearly,  (\ref{Omega*_lj norm}) implies
\bel{F* natural L^1}
\left\|{^\natural}\F a\right\|_{\L^1(\R^\N\setminus\mathfrak{Q})}~\leq~\C~\left\| a\right\|_{\L^1(\R^\N)}.
\eeq
In order to obtain (\ref{H^1 to L^1 EST*}), we repeat an analogue of estimates from (\ref{Est r>1}) to (\ref{Comple est}).  

First, if $\supp a\subset\B_r$ for $r\ge1$, we have
\bel{Est* r>1}
\begin{array}{lr}\ds
\int_{\mathfrak{Q}}\left|{^\natural}\F^* a(x)\right| dx
~\leq~\C~\left\| {^\natural}\F^* a\right\|_{\L^2(\R^\N)}\qquad\hbox{by Schwartz inequality}
\\\\ \ds~~~~~~~~~~~~~~~~~~~~~~~~~
~\leq~\C~\left\| a\right\|_{\L^2(\R^\N)}\qquad \hbox{\small{because ${^\natural}\F$ is bounded on $\L^2(\R^\N)$}}
\\\\ \ds~~~~~~~~~~~~~~~~~~~~~~~~~
~\leq~\C~r^{\N\big[-1+{1\over 2}\big]},\qquad\hbox{\small{ $\|a\|_{\L^2(\R^\N)}\leq r^{-\N+{\N\over 2}}$ }}.
\end{array}
\eeq
Next, for $0<r<1$, the region of influence associated to $\F^*_{\ell}$ is mollified as follows.

Recall that $\left\{ \xi^\nu_j\right\}_\nu$ is a collection of points  equally distributed on  $\mathds{S}^{\N-1}$ with a grid length between $2^{-j/2-1}$ and $2^{-j/2}$. 

Given $\xi^\nu_j\in \mathds{S}^{\N-1}$, consider
\bel{R* xi}
{^*}\hbox{\bf R}^\nu_j~=~\left\{x\in\R^\N~\colon~ \left|x-\left(\nabla_\xi\Phi\right)\left(0,\xi^\nu_j\right)\cdot \xi^\nu_j\right|\leq\c 2^{-j},~~~~\left|x-\left(\nabla_\xi\Phi\right)\left(0,\xi^\nu_j\right)\right|\leq\c 2^{-j/2}\right\}
\eeq
where $\c>0$ is some large constant depending on $\Phi$. 

Define
  \bel{Q*}
\Q_{r\ell}^*~=~\bigcup_{j\colon~2^{-j}\leq r} ~\bigcup_{\nu\colon~\xi_j^\nu\in \mathds{S}^{\N-1}\cap\Lambda_\ell} {^*}\hbox{\bf R}^\nu_j
\eeq
where $\Lambda_\ell$ is the dyadic cone in (\ref{Cone}).

Recall {\bf Remark 4.1}. From (\ref{Q*}), we again have
\bel{Q* natural norm}
\left|\Q^*_{r\ell}\right|~\leq~\C~r\left\{\begin{array}{lr}\ds 2^{-\ell m},\qquad \ell\ge0,
\\ \ds
2^{n\ell},\qquad~~~ \ell\leq0.
\end{array}
\right.
\eeq
By using Schwartz inequality and (\ref{Q* natural norm}), we find
\bel{Local est*}
\begin{array}{lr}\ds
\int_{\Q^*_{r\ell}} \left|{^\natural}\F^*_\ell a(x)\right|dx~\leq~\C\left\|{^\natural}\F^*_\ell a\right\|_{\L^2(\R^\N)}~r^{1\over 2} \left\{\begin{array}{lr}\ds 2^{-m\ell/2},\qquad \ell\ge0,
\\ \ds 2^{n\ell/2},\qquad ~~\ell\leq0.
\end{array}\right.
\end{array}
\eeq
Let ${\rho\over \N}={1\over p}-{1\over 2}$ and $\rho={\N-1\over 2}$.
By applying (\ref{F*_t 2,p result}) in {\bf Proposition One},   we have
\bel{Local est 2*}
\begin{array}{lr}\ds
\int_{\Q^*_{r\ell}} \left|{^\natural}\F^*_\ell a(x)\right|dx~\leq~\C~\left\| a\right\|_{\L^p(\R^\N)} 
r^{1\over 2}2^{\big[{\N\over 2}-\rho-\ve\big]|\ell|}\left\{\begin{array}{lr}\ds 2^{-m\ell/2},\qquad \ell\ge0,
\\ \ds 2^{n\ell/2},\qquad ~~\ell\leq0
\end{array}\right.
\\\\ \ds~~~~~~~~~~~~~~~~~~~~~~~~~~
~\leq~\C~ r^{-\N+{\N\over p}} r^{1\over 2} 2^{\big[{\N\over 2}-\rho-\ve\big]|\ell|}\left\{\begin{array}{lr}\ds 2^{-m\ell/2},\qquad \ell\ge0,
\\ \ds 2^{n\ell/2},\qquad ~~\ell\leq0
\end{array}\right.\qquad\hbox{\small{ $\|a\|_{\L^p(\R^\N)}\leq r^{-\N+{\N\over p}}$ }}
\\\\ \ds~~~~~~~~~~~~~~~~~~~~~~~~~~
~=~\C~ r^{-\N+{\N-1\over 2}+{\N\over 2}} r^{1\over 2} 2^{\big[{1\over 2}-\ve\big]|\ell|}\left\{\begin{array}{lr}\ds 2^{-m\ell/2},\qquad \ell\ge0,
\\ \ds 2^{n\ell/2},\qquad ~~\ell\leq0
\end{array}\right.
\\\\ \ds~~~~~~~~~~~~~~~~~~~~~~~~~~
~=~\C~ 2^{\big[{1\over 2}-\ve\big]|\ell|}\left\{\begin{array}{lr}\ds 2^{-m\ell/2},\qquad \ell\ge0,
\\ \ds 2^{n\ell/2},\qquad ~~\ell\leq0.
\end{array}\right.
\\\\ \ds~~~~~~~~~~~~~~~~~~~~~~~~~~
~\leq~\C~ 2^{-\ve|\ell|},\qquad \ve=\ve(\rho_1,\rho_2)>0.
\end{array}
\eeq
Denote $\Omega^*_{\ell j}$ and ${^\sharp}\Omega^*_j,  {^\flat}\Omega^*_j$ to be the kernel of $\F^*_{\ell j}$ and ${^\sharp}\F^*_j, {^\flat}\F^*_j$ respectively.
\begin{remark}$\Omega^*_{\ell j}$ and $ {^\sharp}\Omega^*_j$, ${^\flat}\Omega^*_j$ satisfy  (\ref{Result Size})-(\ref{Result Size sharp+flat}) with  $\Q_{r\ell}$ replaced by $\Q^*_{r\ell}$.
\end{remark}
This can be verified by carrying out an analogue of  section 4 and 5. 

Lastly, the same argument in (\ref{F natural l rewrite})-(\ref{INT EST sharp+flat}) shows that ${^\natural}\F^*_\ell$ and ${^\sharp}\F^*, {^\flat}\F^*$ satisfy the estimates in 
(\ref{Comple est l}) and (\ref{Comple est}).

\section{Proof of Proposition One}
\setcounter{equation}{0}
We prove the $\L^2$-estimate in (\ref{L^2-result}) by following the lines in {\it 3.1.1}, chapter IX  of Stein \cite{Stein}. 
Consider 
\bel{Sf}
\mathcal{S} f(x)~=~\int_{\R^\N} e^{2\pi\i\Phi(x,\xi)} \omega(x,\xi) f(\xi) d\xi
\eeq 
and 
\bel{S*f}
\mathcal{S}^* f(\xi)~=~\int_{\R^\N} e^{-2\pi\i\Phi(x,\xi)} \bar{\omega}(x,\xi) f(x) dx
\eeq 
where $\omega(x,\xi)$ has a compact support in $x$. Moreover, it satisfies
\bel{sigma x diff}
\left|\partial_x^\alpha\omega(x,\xi)\right|~\leq~\C_\alpha
\eeq
for every multi-index $\alpha$.

Let $\b$  be a small positive constant. We define an {\it narrow cone} as follows: suppose $\xi$ and $\eta$ belong to a same narrow cone and $|\eta|\leq|\xi|$. Write $\eta=\rho\xi+\eta^\dagger$ for $0\leq\rho\leq1$ where $\eta^\dagger$ is perpendicular to $\xi$. Then, we have  $|\eta^\dagger|\leq\b\rho|\xi|$. The value of $\b$ depends on $\Phi$.
Clearly, every $\mathcal{S}$ or $\mathcal{S}^*$ can be written as a finite sum of  partial operators.  Each one of them has  a symbol  supported on an narrow cone.

We borrow the next result from p.397 in the book of Stein \cite{Stein}. Namely,  
\bel{Phi_x est}
\left|\nabla_x \Big[\Phi(x,\xi)-\Phi(x,\eta)\Big]\right|~\ge~\C_\Phi~|\xi-\eta|
\eeq
whenever $\xi$ and $\eta$ belong to a same narrow cone.

A  direct computation shows
\bel{S*Sf}
\begin{array}{cc}\ds
\mathcal{S}^*\mathcal{S} f (\xi)~=~\int_{\R^\N} f(\eta)\K(\xi,\eta)d\eta,
\\\\ \ds
\K(\xi,\eta)~=~\int_{\R^\N} e^{2\pi\i\big[\Phi(x,\eta)-\Phi(x,\xi)\big]} \omega(x,\eta)\bar{\omega}(x,\xi) dx.
\end{array}
\eeq
Because $\omega(x,\xi)$ has a $x$-compact support, $\K(\xi,\eta)$ in (\ref{S*Sf}) is uniformly bounded.  

 By using (\ref{sigma x diff}) and (\ref{Phi_x est}),  
an $N$-fold integration by parts  $w.r.t~x$ shows
\bel{K int by parts}
\begin{array}{lr}\ds
\left|\K(\xi,\eta)\right|~\leq~\C_N~\left({1\over 1+|\xi-\eta|}\right)^N,\qquad N\ge0.
\end{array}
\eeq
From (\ref{S*Sf}), we have
\bel{L^2 est S*S} 
\begin{array}{lr}\ds
\left\|\mathcal{S}^*\mathcal{S} f\right\|_{\L^2(\R^\N)}~=~\left\{\int_{\R^\N}\left|\int_{\R^\N} f(\xi-\eta)\K(\xi,\xi-\eta)d\eta\right|^2d\xi\right\}^{1\over 2}\qquad\hbox{\small{$\eta\mt\xi-\eta$}}
\\\\ \ds~~~~~~~~~~~~~~~~~~~~~
~\leq~\C~\int_{\R^\N} \left\{\int_{\R^\N} \left|f(\xi-\eta)\right|^2 \Big|\K(\xi,\xi-\eta)\Big|^2 d\xi \right\}^{1\over 2 }d\eta~~~\hbox{\small{by Minkowski integral inequality}}
\\\\ \ds~~~~~~~~~~~~~~~~~~~~~
~\leq~\C~\left\| f\right\|_{\L^2(\R^\N)}\int_{\R^\N}  \left({1\over 1+|\eta|}\right)^{\N+1}d\eta\qquad \hbox{\small{by (\ref{K int by parts}) with $N=\N+1$}}.
\end{array}
\eeq
Let $\F_o$ and ${^\natural}\F, {^\sharp}\F, {^\flat}\F$  defined in (\ref{F_o}) and (\ref{F natural})-(\ref{F flat}). Each one of them is a composition of Fourier transform and the operator $\mathcal{S}$ defined in (\ref{Sf}). Moreover, $\omega(x,\xi)$ equals
\bel{symbols}
\begin{array}{cc}\ds
\phi(\xi)\sigma(x,\xi),\qquad \sigma(x,\xi)\sum_{j>0}\phi_j(\xi)\sum_{-j/2<\ell<j/2}\delta_\ell(\xi)
\\\\ \ds
\sigma(x,\xi)\sum_{j>0}\phi_j(\xi) \sum_{\ell\ge j/2}\delta_\ell(\xi),\qquad \sigma(x,\xi)\sum_{j>0}\phi_j(\xi) \sum_{\ell\leq -j/2}\delta_\ell(\xi)
\end{array}
\eeq
respectively where $\delta_\ell$ and $\phi, \phi_j$ are defined in (\ref{delta_t}) and (\ref{phi_j}). 

Every symbol function in (\ref{symbols}) satisfies (\ref{sigma x diff}) provided that $\sigma\in\S^0$. Plancherel theorem together with (\ref{L^2 est S*S}) imply (\ref{L^2-result}).

Next, we write $\xi=(\tau,\lambda)\in\R^n\times\R^m$ and $x=(z,w)\in\R^n\times\R^m$, $y=(u,v)\in\R^n\times\R^m$. 
Consider
\bel{T}
\begin{array}{cc}\ds
  \T f(x)~=~\int_{\R^\N}f(y)\left\{\int_{\R^\N} e^{2\pi\i (x-y)\cdot\xi} |\tau|^{-\rho_1}|\lambda|^{-\rho_2} d\xi\right\}dy
 \end{array}
\eeq
where $\rho_1>\left[{n-1\over \N-1}\right]\rho$, $\rho_2>\left[{m-1\over \N-1}\right]\rho$ and $\rho_1+\rho_2=\rho>0$. 

Recall {\bf Remark 2.3}. 
For $0<\rho\leq(\N-1)/2$, we must have $0<\rho_1<{n\over 2}$ and $0<\rho_2<{m\over 2}$.

The kernel of $\T$ defined in (\ref{T}) is a distribution equal to 
\bel{kernel T}
{\pi^{\rho}\Gamma(n-\rho_1)\Gamma(m-\rho_2)\over \pi^{\N-\rho}\Gamma(\rho_1)\Gamma(\rho_2)} \left({1\over |u|}\right)^{n-\rho_1}\left({1\over |v|}\right)^{m-\rho_2},\qquad |u|>0,~~|v|>0.
\eeq
Let ${\rho_1\over n}={1\over p_1}-{1\over 2}$ and  ${\rho_2\over m}={1\over p_2}-{1\over 2}$. By applying Hardy-Littlewood-Sobolev inequality \cite{Hardy-Littlewood}-\cite{Sobolev} in $\R^n$ and $\R^m$ respectively, 
we find
\bel{T EST mixnorm}
\begin{array}{lr}\ds
\left\| \T f\right\|_{\L^2(\R^\N)}~\leq~\C_{\rho_1~\rho_2} \left\{\iint_{\R^n\times\R^m} \left\{\iint_{\R^n\times\R^m} |f(z-u,w-v)| \left({1\over |u|}\right)^{n-\rho_1}\left({1\over |v|}\right)^{m-\rho_2} dudv\right\}^2 dzdw\right\}^{1\over 2}
\\\\ \ds~~~~~~~~~~~~~~~~~~~
~\leq~\C_{\rho_1~\rho_2} \left\{\int_{\R^m} \left\{\int_{\R^n}\left\{\int_{\R^m} |f(z,w-v)| \left({1\over |v|}\right)^{m-\rho_2} dv\right\}^{p_1} dz\right\}^{2\over p_1}dw\right\}^{1\over 2}
\\\\ \ds~~~~~~~~~~~~~~~~~~~
~\leq~\C_{\rho_1~\rho_2} \left\{\int_{\R^n} \left\{\int_{\R^m}\left\{\int_{\R^m} |f(z,w-v)| \left({1\over |v|}\right)^{m-\rho_2} dv\right\}^2 dw\right\}^{p_1\over 2}dz\right\}^{1\over p_1}
\\ \ds~~~~~~~~~~~~~~~~~~~~~~~~~~~~~~~~~~~~~~~~~~~~~~~~~~~~~~~~~~~~~~~~~~~~~~~~~
\hbox{\small{by Minkowski integral inequality}}
\\ \ds~~~~~~~~~~~~~~~~~~~
~\leq~\C_{\rho_1~\rho_2} \left\{\int_{\R^n} \left\{\int_{\R^m} |f(z,w)|^{p_2}  dw\right\}^{p_1\over p_2}dz\right\}^{1\over p_1}.
\end{array}
\eeq
Let $\F_\ell$ defined in (\ref{Delta_lj F})  and ${^\sharp}\F, {^\flat}\F$ defined in (\ref{F sharp})-(\ref{F flat}). We write 
\bel{F_l rewrite}
\begin{array}{lr}\ds
\F_\ell f(x)~=~\int_{\R^\N} e^{2\pi\i\Phi(x,\xi)}\left[\delta_\ell(\xi)\sigma(x,\xi)|\tau|^{\rho_1}|\lambda|^{\rho_2} \sum_{j>2|\ell|}\phi_j(\xi) \right]\Hat{\T f}(\xi)d\xi
\end{array}
\eeq
and
\bel{F sharp/flat rewrite}
\begin{array}{lr}\ds
{^\sharp}\F f(x)~=~\int_{\R^\N} e^{2\pi\i\Phi(x,\xi)}\left[\sigma(x,\xi)|\tau|^{\rho_1}|\lambda|^{\rho_2}\sum_{j>0}\phi_j(\xi) \sum_{\ell\ge j/2}\delta_\ell(\xi)\right] \Hat{\T f}(\xi)d\xi,
\\\\ \ds
{^\flat}\F_j f(x)~=~\int_{\R^\N} e^{2\pi\i\Phi(x,\xi)}\left[\sigma(x,\xi)|\tau|^{\rho_1}|\lambda|^{\rho_2}\sum_{j>0}\phi_j(\xi) \sum_{\ell\leq -j/2}\delta_\ell(\xi)\right] \Hat{\T f}(\xi)d\xi.
\end{array}
\eeq
Suppose $\sigma\in\S^{-\rho}$ satisfying the differential inequality in (\ref{Class})-(\ref{rho_1, rho_2}). 
Every symbol function $\left[\cdots\right]$ in (\ref{F_l rewrite})-(\ref{F sharp/flat rewrite}) 
 satisfies  (\ref{sigma x diff}) consequently.
Therefore, we obtain  (\ref{F_t 2,p result})-(\ref{formula rho_1, rho_2}) by using the regarding $\L^2$-estimate together with (\ref{T EST mixnorm}).

Finally, we prove (\ref{F*_t 2,p result}). Denote
\bel{Delta_l}
\Delta_\ell (\xi)~=~\delta_\ell(\xi)\sum_{j>2|\ell|}\phi_j(\xi),\qquad \ell\in\Z. 
\eeq
Consider
\bel{F*_l}
\begin{array}{cc}\ds
{^\natural}\F^*_\ell f(x)~=~\int_{\R^\N} e^{2\pi\i x\cdot\xi} {^\natural}\mathcal{S}^*_\ell f(\xi)d\xi,
\\\\ \ds
{^\natural}\mathcal{S}^*_\ell f(\xi)~=~\int_{\R^\N} e^{-2\pi\i\Phi(x,\xi)} \bar{\Delta}_\ell(\xi)\bar{\sigma}(x,\xi) f(x) dx.
\end{array}
\eeq
Let $\sigma\in\S^{-\rho}, 0<\rho\leq(\N-1)/2$. Each ${^\natural}\mathcal{S}^*_\ell$ is bounded on $\L^2(\R^\N)$. Indeed,  $\bar{\Delta}_\ell(\xi)\bar{\sigma}(x,\xi)$ again satisfies  (\ref{sigma x diff}). Because of Plancherel theorem, it is suffice to show
\bel{S*_l regularity}
\left\| {^\natural}\mathcal{S}^*_\ell f\right\|_{\L^2(\R^\N)}~\leq~\C~2^{\big[{\N\over 2}-\rho-\ve\big]|\ell|}~\left\|f\right\|_{\L^p(\R^\N)},\qquad {\rho\over \N}~=~{1\over p}-{1\over 2}.
\eeq
By symmetry, we assert $\ell\ge0$ only.
H\"{o}lder inequality implies
\bel{S^*_t L^2-norm}
\begin{array}{lr}\ds
\left\| {^\natural}\mathcal{S}^*_\ell f\right\|_{\L^2(\R^\N)}^2~\leq~\left\|{^\natural}\mathcal{S}_\ell{^\natural}\mathcal{S}^*_\ell f\right\|_{\L^{p\over p-1}(\R^\N)}\left\| f\right\|_{\L^p(\R^\N)}.
\end{array}
\eeq
We prove (\ref{S*_l regularity}) by showing
\bel{SS*_l regularity}
\begin{array}{cc}\ds
\left\| {^\natural}\mathcal{S}_\ell{^\natural}\mathcal{S}^*_\ell f\right\|_{\L^{p\over p-1}(\R^\N)}~\leq~\C ~2^{\big[\N-2\rho-2\ve\big]\ell}~\left\| f\right\|_{\L^p(\R^\N)},
\qquad
{2\rho\over \N}~=~{1\over p}-{p-1\over p}.
\end{array}
\eeq
From direct computation, we find that kernel of ${^\natural}\mathcal{S}_\ell{^\natural}\mathcal{S}^*_\ell$ 
 is a distribution equal to
\bel{V_l}
\V_\ell(x,y)~=~\int_{\R^\N} e^{2\pi\i\big[\Phi(x,\xi)-\Phi(y,\xi)\big]} \sigma(x,\xi)\Delta_\ell(\xi)\bar{\sigma}(y,\xi)\bar{\Delta}_\ell(\xi)d\xi
\eeq
whenever it is well defined at $(x,y)\in\R^\N\times\R^\N$.

Let $\phi_k, k\in\Z$ defined in (\ref{phi_j}). Recall that $\Phi$ is homogeneous of degree $1$ in $\xi$. Consider
\bel{V_lk}
\begin{array}{lr}\ds
\V_{\ell k}(x,y)~=~\int_{\R^\N} e^{2\pi\i\big[\Phi(x,\xi)-\Phi(y,\xi)\big]} \sigma(x,\xi)\Delta_\ell(\xi)\bar{\sigma}(y,\xi)\bar{\Delta}_\ell(\xi)\phi_k(\xi)d\xi
\\\\ \ds~~~~~~~~~~~~~~
~=~2^{k\N} \int_{\R^\N} e^{2\pi\i2^k\big[\Phi(x,\xi)-\Phi(y,\xi)\big]} \sigma\left(x,2^k\xi\right)\Delta_\ell\left(2^k\xi\right)\bar{\sigma}\left(y,2^k\xi\right)\bar{\Delta}_\ell\left(2^k\xi\right)\phi_o(\xi)d\xi
\\ \ds~~~~~~~~~~~~~~~~~~~~~~~~~~~~~~~~~~~~~~~~~~~~~~~~~~~~~~~~~~~~~~~~~~~~~~~~~~~~~~~~~~~~~~~~~~~~~~~~~~~~~~
\hbox{\small{$\xi\mt 2^k \xi$,~~~~~$k\in\Z$}}.
\end{array}
\eeq
Note that $\delta_\ell(2^k\xi)=\delta_\ell(\xi)$ as defined in (\ref{delta_t}).
From (\ref{Delta_l}), we find
\bel{Delta_l k-dila}
\begin{array}{lr}\ds
\Delta_\ell\left(2^k\xi\right)\bar{\Delta}_\ell\left(2^k\xi\right)\phi_o(\xi)
~=~\delta_\ell(\xi) \bar{\delta}_\ell(\xi)\phi_o(\xi)\left[\sum_{j>2|\ell|\colon j=k-1,k,k+1}\phi_j\left(2^k\xi\right)\right]
\left[\sum_{j>2|\ell|\colon j=k-1,k,k+1}\bar{\phi}_j\left(2^k\xi\right)\right]. 
\end{array}
\eeq
Write $\Phi_{x\xi}(x,\xi)$ for the $\N\times \N$-matrix $\left[{\p^2\Phi\over\p x\p\xi}\right](x,\xi)$.
We have
\bel{nabla Phi}
\nabla_\xi \Big[\Phi(x,\xi)-\Phi(y,\xi)\Big]~=~\Phi_{x\xi}(x,\xi)(x-y)+\O(1)|x-y|^2.
\eeq
\begin{remark}
Momentarily, we assume that $\sigma\in\S^\rho$ has a sufficiently small  $x$-support.
\end{remark}
From (\ref{nabla Phi}), we have
\bel{nabla Phi size}
\left|\nabla_\xi \Big[\Phi(x,\xi)-\Phi(y,\xi)\Big]\right|~\ge~\C~|x-y|
\eeq
for $|x-y|$ sufficiently small as $\Phi(x,\xi)$ satisfying the non-degeneracy condition in (\ref{nondegeneracy}).

Suppose $\sigma\in\S^{-\rho}, 0<\rho\leq(\N-1)/2$.   
Recall {\bf Remark 2.3}. 
We have $0<\rho_1<{n\over 2}$ and $0<\rho_2<{m\over 2}$ as an necessity.
For $\xi\in\supp\delta_\ell\left(\xi\right) \bar{\delta}_\ell\left(\xi\right)\phi_o\left(\xi\right)$, we find $|\tau|\sim1$, $|\lambda|\sim2^{-\ell}$ and
\bel{symbols size}
\begin{array}{lr}\ds
\left| \sigma\left(x,2^k\xi\right)\bar{\sigma}\left(y,2^k\xi\right)\right|~\leq~\left({1\over 1+2^k }\right)^{2\rho_1}  \left({1\over 1+2^{k-\ell} }\right)^{2\rho_2}
\\\\ \ds~~~~~~~~~~~~~~~~~~~~~~~~~~~~~~~~~~
~\leq~\C~2^{-2\rho k} 2^{2 \rho_2\ell}
~=~\C~2^{-2\rho k} 2^{m\ell}2^{-2\ve\ell}
\end{array}
\eeq
where $2\ve=m-2\rho_2>0$. 
On the other hand, $\delta_\ell$ satisfies the differential inequality in (\ref{delta_l Diff Ineq}). 
From direct computation, we find
\bel{Diff Ineq flat j}
\begin{array}{lr}\ds
\left| \nabla_\xi^N
\sigma\left(x,2^k\xi\right)\bar{\sigma}\left(y,2^k\xi\right)\Delta_\ell\left(2^k\xi\right)\bar{\Delta}_\ell\left(2^k\xi\right)\phi_o(\xi)\right|
~\leq~\C_N ~2^{-2\rho k} 2^{m\ell}2^{-2\ve\ell} 2^{\ell N}
\end{array}
\eeq
for every $N\ge0$.

Let $\V_{\ell j}(x,y)$ defined in (\ref{V_lk}). By using (\ref{nabla Phi size}) and (\ref{Diff Ineq flat j}),  an $N$-fold  integration by parts $w.r.t~ \xi$ shows
\bel{V_lk by parts}
\begin{array}{lr}\ds
\left|\V_{\ell j}(x,y)\right|
~\leq~\C~ \Big[ 2^{-2\rho k} 2^{m\ell} 2^{-2\ve\ell}\Big] \Big[2^{k\N}2^{-\ell m}\Big] ~\left[2^{k-\ell} |x-y|\right]^{-N}
\\\\ \ds~~~~~~~~~~~~~~~~~
~=~ \C~ 2^{-2\ve\ell}2^{\big[\N-2\rho\big]k} \left[2^{k-\ell} |x-y|\right]^{-N}
\\\\ \ds~~~~~~~~~~~~~~~~~
~=~\C~2^{-2\ve\ell}2^{\big[\N-2\rho\big]\ell}~2^{\big[\N-2\rho\big](k-\ell)} \left[2^{k-\ell} |x-y|\right]^{-N}.\end{array}
\eeq
Choose 
\bel{N,M lambda}
\begin{array}{cc}\ds
N=0~~~\hbox{if}~~~ |x-y|\leq2^{-k+\ell}\qquad \hbox{or}\qquad N=n~~~\hbox{if}~~~ |x-y|>2^{-k+\ell}.
\end{array}
\eeq
Recall $\V_\ell$ defined in (\ref{V_l}) which is the kernel of ${^\natural}\mathcal{S}_\ell{^\natural}\mathcal{S}^*_\ell$.  By using (\ref{V_lk by parts})-(\ref{N,M lambda}), we find
\bel{K_l Size}
\begin{array}{lr}\ds
\left|\V_\ell(x,y)\right|
~\leq~\C~2^{\big[\N-2\rho-2\ve\big]\ell} \left({1\over |x-y|}\right)^{\N-2\rho},\qquad\hbox{\small{$|x-y|>0$}}.
\end{array}
\eeq
By applying Hardy-Littlewood-Sobolev inequality \cite{Hardy-Littlewood}-\cite{Sobolev}, we obtain (\ref{SS*_l regularity}).

Lastly, the extra assumption on $\supp\sigma$ made in {\bf Remark 7.1} can be removed because $\sigma(x,\xi)$ has a compact support in $x$.

 {\small Department of Mathematics, Shantou University}
  \\
 {\small email: cqtan@stu.edu.cn}

 {\small Department of Mathematics, Westlake University} \\
 {\small email: wangzipeng@westlake.edu.cn}

\end{document}